\documentclass[11pt,reqno,twoside]{article}
\usepackage{amsmath,amsfonts,amssymb}
\usepackage{mathrsfs}
\usepackage{amsthm}
\usepackage{enumitem}
\usepackage[utf8]{inputenc}
\usepackage{fontenc}
\usepackage{hyperref}
\hypersetup{  
   bookmarks=true,
   backref=true,
   pagebackref=false,
   colorlinks=true,
   linkcolor=blue,
   citecolor=red,
   urlcolor=blue
}

\usepackage{subfig,float}
\usepackage[english]{babel}
\usepackage{amsbsy,amscd,fancyhdr,graphicx,psfrag,fancybox,indentfirst,color}
\usepackage{graphics,epsfig}
\usepackage{mysty_arxiv}
\usepackage{fullpage}
\usepackage{graphicx}

\newcommand*\samethanks[1][\value{footnote}]{\footnotemark[#1]}

\newcommand{\hl}[1]{{\textcolor{black}{#1}}}

\title{\textbf {Nonlocal $p$-Laplacian Variational problems on graphs}}         
\date{}
\author{Yosra Hafiene\thanks{Normandie Univ, ENSICAEN, UNICAEN, CNRS, GREYC, France.} \and Jalal M. Fadili\samethanks \and Abderrahim Elmoataz\samethanks[1]
}
\begin{document}

\maketitle
\begin{abstract}
In this paper, we study a nonlocal variational problem which consists of minimizing in $L^2$ the sum of a quadratic data fidelity and a regularization term corresponding to the $L^p$-norm of the nonlocal gradient. In particular, we study convergence of the numerical solution to a discrete version of this nonlocal variational problem to the unique solution of the continuum one. To do so, we derive an error bound and highlight the role of the initial data and the kernel governing the nonlocal interactions. When applied to variational problem on graphs, this error bound allows us to show the consistency of the discretized variational problem as the number of vertices goes to infinity. More precisely, for networks in convergent graph sequences (simple and weighted deterministic dense graphs as well as random inhomogeneous graphs), we prove convergence and provide rate of convergence of solutions for the discrete models to the solution of the continuum problem as the number of vertices grows.  
%In this paper we consider the regularization problem associated to the nonlocal $p$-Laplacian operator. We study the consistency of its numerical approximations
\end{abstract}


\section{Introduction}
\subsection{Problem statement}
%The clear image is defined by the following variational problem 
\hl{Let $\O \subset \R$ be a bounded domain, and without loss of generality we take $\O = [0,1]$, and the kernel $K$ is a symmetric, nonnegative, measurable and bounded function on $\O^2$. Denote the nonlocal gradient operator 
\[
\nabla_K: u \in L^2(\O) \mapsto \pa{K(x,y)^{1/p}(u(y)-u(x))}_{(x,y) \in \O^2} \in L^2(\O^2). 
\]
We study the following variational problem for $ p \in [1, + \infty[$
\begin{equation}\tag{\textrm{$\mathcal{VP}^{\lambda, p}$}}
\min_{u \in L^2(\O)} \acc{ \Elam (u, g, K) \eqdef    \frac{1}{2\lambda} \norm{u-g}_{L^2(\O)}^2 + \reg(u, K)} .
\label{def : mainpb}
\end{equation}
The nonlocal regularizer $\reg$ is defined as
\begin{equation}
\reg(u, K) \eqdef 
\begin{cases}
\frac{1}{2p}\int_{\O^2} \abs{\nabla_K u (x,y)}^p dx dy & \text{if } |\cdot|^p \circ \nabla_K u \in L^1(\O^2), \\
+\infty & \text{otherwise} .
\end{cases}
\label{def : flow}
\end{equation}}
%We deal with the following problem 
%\[
%E : L^2(\O) \times L^2(\O) \times L^{\infty}(\O) \longrightarrow \R 
%\]
%\[
%\hspace{6cm} (u, g, K) \mapsto \frac{1}{2\lambda} \norm{u-g}_{L^2(\O)} + \int_{\O^2} K(x,y) \abs{u(y)-u(x)}^p dx dy 
%\]
%\begin{equation}
%E(u, g, K) = 
Here $\lambda$ is a positive regularization parameter. The chief goal of this paper is to study numerical approximations of the nonlocal variational problem~\eqref{def : mainpb}, which in turn, will allow us to establish consistency estimates of the discrete counterpart of this problem on graphs. 
%Let us note that the gradient of the functional $\reg (u, K)$ is given by the nonlocal $p$-Laplacian operator 
%\begin{equation}
% \A(u(x)) =  - \displaystyle{\int}_{\O} \K(x,y) \abs{u(y) - u(x) }^{p-2} (u(y) - u(x) ) dy.
% \label{plaplaceoperator}
%\end{equation}$

%ici j'ai rajouté un paragraphe d'intro pour les motivations
\hl{Discretization of continuum models based on nonlocal regularization such as in \eqref{def : mainpb} has proven very effective in various tasks in signal/image/data processing, machine learning and computer vision. Such models have the advantage of better mathematical modeling, connections with physics and better geometrical approximations.}

\hl{In many real-world problems where the data in practice is discrete, graphs constitute a natural structure suited to their representation. Each vertex of the graph corresponds to a datum, and the edges encode the pairwise relationships or similarities among the data. For the particular case of images, pixels (represented by nodes) have a specific organization expressed by their spatial connectivity. Therefore, a typical graph used to represent images is a grid graph. For the case of unorganized data such as point clouds, a graph can also be built by modeling neighborhood relationships between the data elements. For these practical reasons, recently, there has been a surge of interest in adapting and solving
nonlocal variational problems such as \eqref{def : mainpb} on data which is represented by arbitrary graphs and networks. Using this framework, problems are directly expressed in a discrete
setting. This way to proceed encompasses local and nonlocal methods in the same framework by using
appropriate graph topologies and edge weights depending on the data structure and the task to
be performed.}

\hl{Thus, handling such data necessitates a discrete counterpart of~\eqref{def : mainpb}. One can intuitively propose the following problem
\begin{equation}\tag{\textrm{$\mathcal{VP}_n^{\lambda, p}$}}
\min_{u_n \in \R^n}\acc{ \Elamn \eqdef  \frac{1}{2 \lambda n } \norm{u_n - g_n}_2^2 + \regd(u_n,K_n) } ,
\label{def : maindiscpb}
\end{equation}
where
\begin{equation}
\regd(u_n, K_n) \eqdef \frac{1}{2n^2p} \sum\limits_{i,j = 1}^{n} K_{nij} \abs{u_{nj} - u_{ni}}^p .
\label{def : flowd}
\end{equation}
Since the discrete nonlocal problem~\eqref{def : mainpb} was only intuitively inspired as a numerical approximation to \eqref{def : mainpb}, several legitimate questions then arise:
\begin{enumerate}[label=(Q.\arabic*)]
\item Is there any continuum limit to the (unique) minimizer $\unast$ of \eqref{def : maindiscpb} as $n \to +\infty$ ? If yes, in what sense ? The graphs considered in \eqref{def : maindiscpb} are rather general and it is not obvious at first glance what a continuum model should be.
\item What is the rate of convergence to this limit and what is its relation to the (unique)
minimizer $\uast$ of \eqref{def : mainpb} ? 
\item What are the parameters involved in this convergence rate (in particular the interplay between $n$, $p$ and properties of the graph), and what is their influence in the corresponding rate ?
\item Can this continuum limit help us get better insight into discrete models/algorithms and how they may be improved ?
\end{enumerate}
It is our primary motivation to answer these questions rigorously. Our perspective stands at the interface of numerical analysis and data processing.}

\subsection{Contributions}
In this work we focus on studying the consistency of~\eqref{def : mainpb} in which we investigate functionals with a nonlocal regularization term corresponding to the $p$-Laplacian operator. We first give a general error estimate in $L^2(\O)$ controlling the error between the continuum extension of the numerical solution $\unast$ to the discrete variational problem~\eqref{def : maindiscpb} and its continuum analogue $\uast$ of~\eqref{def : mainpb}. The dependence of the bound on the error induced by discretizing the kernel $K$ and the initial data $g$ is made explicit. Under very mild conditions on $K$ and $g$, typically belonging to a large class of Lipschitz functional spaces (see Section~\ref{subsec: Lipspaces} for details on these spaces), convergence rates can be exhibited \hl{by plugging Lemma~\ref{lem:lipgu} and Lemma~\ref{lem:spaceapprox} into the bounds of Theorem~\ref{mainresult}. It is worth pointint out that at this stage, though we have focused on the 1D case, these results hold true if $\O$ is any bounded domain in $\R^d$. This would allow for instance the application of our results to a variety of situations, typically image or point cloud processing (see Section~\ref{sec:numerical} for a detailed example). Here the nodes/vertices of the graph are the pixels/points locations, and the edge weights model the neighborhood relationships. The neighbourhood can be thought of as spatial or also include value similarity such as in patch-based methods.}

Secondly, we apply these results, using the theory graph limits (for instance graphons), to dynamical networks on simple and weighted dense graphs to show that the approximation of minimizers of the discrete problems on simple and weighted graph sequences converge to those of the continuum problem. This sets the question that solving a discrete variational problem on graphs has indeed a continuum limit. 
%We give also a rate of convergence estimate. 
Under very mild conditions on $K$ and $g$, typically belonging to Lipschitz functional spaces, precise convergence rates can be exhibited. These functional spaces allow one to cover a large class of graphs (through $K$) and initial data $g$, including those functions of bounded variation.
For simple graph sequences, we also show how the accuracy of the approximation depends on the regularity of the boundary of the support of the graph limit.
 %which in particular contain functions of bounded variation. 

Finally, building upon these error estimates, we study networks on random inhomogeneous graphs. We combine them with sharp deviation inequalities to establish nonasymptotic convergence claims and give the rate of convergence of the discrete solution to its continuum limit with high probability under the same assumptions on the kernel $K$ and the initial data $g$. 
%and show the influence of the value of p
%As random graphs arise in many applications 

\subsection{Relation to prior work}
%\subsection{Discussion of related works}
\paragraph{Graph-based regularization in machine learning} 
\hl{Semi-supervised learning with a weighted graph to capture the geometry of the unlabelled data and graph Laplacian-based regularization is now popular; see \cite{vonLuxburg07} for a review. Pointwise and/or spectral convergence of graph Laplacians were studied by several authors including \cite{Belkin07,Hein05,Hein06,Coifman06,Singer06,vonLuxburg08,Ting10,Singer17}.}

\hl{Graph-based $p$-Laplacian regularization has also found applications in semi-supervised learning such as clustering; see e.g. \cite{Zhou05,elmoataz08,Buhler09}. For instance, the authors of~\cite{wainwright16} obtained iterated pointwise convergence of {\textit{rescaled}} graph $p$-Laplacian energies to the continuum (local) $p$-Laplacian as the fraction of labelled to unlabelled points is vanishingly small. The authors in~\cite{SlepcevGarcia2016} studied the consistency of {\textit{rescaled}} total variation minimization on random point clouds in $\R^d$ for a clustering application. They considered the total variation on graphs with a radially symmetric and rescaled kernel $K(x,y)=\varepsilon^{-d}J(|x-y|/\varepsilon)$, $\varepsilon > 0$. This corresponds to an instance of $\regd$ for $d=1$ and $p=1$. Under some assumptions on $J$, and for an appropriate scaling of $\varepsilon$ with respect to $n$, which makes the method become localised in the large data limit, they proved that the discrete total variation on graphs $\Gamma$-converges in an appropriate topology, as $n \to \infty$, to weighted local total variation, where the weight function is the density of the point cloud distribution. Motivated by the work of~\cite{wainwright16}, the authors of~\cite{SlepcevThrope2017} studied consistency of the graph $p$-Laplacian for semi-supervised learning in $\R^d$. They considered both constrained and penalized minimization of $\regd$ with a radially symmetric and rescaled kernel as explained before. They uncovered regimes of $p$ and ranges on the scaling of $\varepsilon$ with respect to $n$ for the asymptotic consistency (in the sense of $\Gamma$-convergence) to hold. Continuing along the lines of~\cite{wainwright16}, the work of \cite{Calder2017} studies the consistency of Lipschitz semi-supervised learning (i.e., $p \to \infty$) on graphs in the same asymptotic limit. This work proves that Lipschitz learning is well-posed and showed that the learned functions converge to the solution of an $\infty$-Laplace type equation, depending on the choice of weights in the graph.}

\hl{The Ginzburg-Landau functional has been adapted in \cite{Bertozzi12} to weighted graphs in an application to machine learning and data clustering. In \cite{vanGennip12}, the authors study $\Gamma$-convergence of the graph based Ginzburg-Landau functional, both the limit for zero diffusive interface parameter or when the number of graph vertices increases.}

\hl{However, all the above results are asymptotic and do not provide any error estimates for finite $n$.}

\paragraph{Nonlocal and graph-based regularization in imaging} 
Several edge-aware filtering schemes have been proposed in the literature~\cite{Yaroslavsky85,Smith97,Tomasi98,Kimmel07}. The nonlocal means filter~\cite{buadesnlmeans05} averages pixels that can be arbitrary far away, using a similarity measure based on distance between patches. As shown in \cite{Kindermann2006,Peyre08}, these filters can also be interpreted within the variational framework with nonlocal regularization functionals. They correspond to one step of gradient descent on~\eqref{def : maindiscpb} with $p=2$, where $K_{nij}=J(|x_i-x_j|)$ is computed from the noisy  input image $g$ using either a distance between the pixels $x_i$ and $x_j$~\cite{Yaroslavsky85,Tomasi98,Kimmel07} or a distance between the patches around $x_i$ and $x_j$ \cite{buadesnlmeans05,Szlam08}. This nonlocal variational denoising can be related to sparsity in an adapted basis of eigenvector of the nonlocal diffusion operator~\cite{Coifman05,Szlam08,Peyre08}. This nonlocal variational framework was also extended to handle several linear inverse problems~\cite{Kindermann2006,Gilboa07,Buades07nonlocal,Gilboa08,elmoataz08,elmoataz15}. In~\cite{Peyre11,Facciolo09,Yang13}, the authors proposed a variational framework with nonlocal regularizers on graphs to solve linear inverse problems in imaging where both the image to recover and the graph structure are inferred.

\paragraph{Consistency of the ROF model} 
For local variational problems, the only work on consistency with error bounds that we are aware of is the one of~\cite{Wang2011} who studied the numerical approximation of the Rudin-Osher-Fatemi (ROF) model, which amounts to minimizing in $L^2(\O^2)$ the well-known energy functional 
\[
E(v) \eqdef \frac{1}{2\lambda} \norm{u-g}_{L^2(\O^2)}^2 +  \norm{v}_{\TVnorm(\O^2)}, 
\]
where $g \in L^2(\O^2)$, and $\norm{\cdot}_{\TVnorm(\O^2)}$ denotes the total variation seminorm. They bound the difference between the continuum solution and the solutions to various finite-difference approximations (including the upwind scheme) to this model. They gave an error estimate in $L^2(\O^2)$ of the difference between these two solutions and showed that it scales as $n^{-\tfrac{s}{2(s+1)}}$, where $s \in ]0,1]$ is the smoothness parameter of the Lipschitz space containing $g$.
  
However, to the best of our knowledge, there is no such error bounds in the nonlocal variational setting. In particular, the problem of the continuum limit and consistency of~\eqref{def : maindiscpb} with error estimates is still open in the literature. It is our aim in this work to rigorously settle this question.

\subsection{Paper organisation}
The rest of this paper is organized as follows. Section~\ref{sec:backmaterialvariational} collects some notations and preliminaries that we will need in our exposition. In Section~\ref{sec:exisvariational} we briefly discuss well-posedness of problems~\eqref{def : mainpb} and~\eqref{def : maindiscpb} and recall some properties of the corresponding minimizers. Section~\ref{sec:errorestimatevariational} is devoted to the main result of the paper (Theorem~\ref{mainresult}) in which we give a bound on the $L^2$-norm of the difference between the unique minimizers of~\eqref{def : mainpb} and~\eqref{def : maindiscpb}. In this section, we also state a key regularity result on the minimizer $\uast$ of~\eqref{def : mainpb}. This result is then used to study networks on deterministic dense graph sequences in Section~\ref{sec:deteministicapplivariational}. First we deal with networks in simple graphs, and show in Corollary~\ref{cor:simplegraphs} the influence of the regularity of the boundary of the support of the graphon on the convergence rate. Secondly, in Section~\ref{sec:weightedvariational} we study networks on weighted graphs. Section~\ref{sec:randomapplivariational} deals with networks on random inhomogeneous graphs. We quantify the rate of convergence with high probability. Numerical results are finally reported in Section~\ref{sec:numerical} to illustrate our theoretical findings.

%%%%%%%%%%%%%%%%%%%%%%%%%%%%%%%%%%%%%%%%%%%%%%
\section{Notations and preliminaries}
\label{sec:backmaterialvariational}
To provide a self-contained exposition, we will recall two key frameworks our work relies on. The first is the limit graph theory which is the notion of convergence for graph sequences developed for the analysis of networks on graphs. The second is that of Lipschitz spaces that will be instrumental to quantify the rate of convergence in our error bounds. %In our review, we keep our explanations as short as possible and we omit technical details in order not to disturb the flow of the main ideas.

\subsection{Projector and injector}
Let $n \in \N^*$, and divide $\O$ into $n$ intervals
\[
\O_1^{(n)} =  \left [0, \frac{1}{n} \right [ , \O_2^{(n)} =  \left [\frac{1}{n}, \frac{2}{n} \right [,\ldots,\O_j^{(n)} =  \left [\frac{j-1}{n}, \frac{j}{n} \right [,\ldots,\O_n^{(n)} =  \left [ \frac{n-1}{n} , 1\right [,
\]
and let $\mathcal{Q}_n$ denote the partition of $\O$, $\mathcal{Q}_n = \{ \O_i^{(n)}, i \in [n] \eqdef \acc{1,\cdots,n}\}$. Denote $\O_{ij}^{(n)} \eqdef \O_{i}^{(n)} \times \O_{j}^{(n)}$. Without loss of generality, we assume that the points are equispaced so that $|\O_i^{(n)}|=1/n$, where $|\O_i^{(n)}|$ is the measure of $\O_i^{(n)}$. The discussion can be easily extended to non-equispaced points by appropriate normalization; see~Section~\ref{sec:randomapplivariational}.

We also consider the operator $\projn: L^1(\O) \to \R^n$
\[
(\projn v)_{i} \eqdef  \frac{1}{|\O_i^{(n)}|} \int_{\O_i^{(n)}} v(x) dx.
\]
This operator can be also seen as a piecewise constant projector of $u$ on the space of discrete functions. For simplicity, and with a slight abuse of notation, we keep the same notation for the projector $\projn: L^1(\O^2) \to \R^{n \times n}$.

We assume that the discrete initial data $g_n$ and the discrete kernel $K_n$ are constructed as  
\begin{equation}
g_n = \projn g \eqdef \pa{g_{n1}, \cdots, g_{nn}}^\top \qandq  \Kn = \projn K \eqdef (K_{nij})_{1\leq i,j \leq n} ,
\label{projectorgeneral}
\end{equation}
where
\begin{equation}
g_{ni} = (\projn g)_{i} = \frac{1}{|\O_i^{(n)}|} \int_{\O_i^{(n)}} g(x) dx  \qandq K_{nij} = (\projn K)_{ij} = \frac{1}{|\O_{ij}^{(n)}|} \int_{\O_{ij}^{(n)}} K(x,y) dxdy .
\label{eq:averageg}
\end{equation}
% We fix $n \in \N^*$, divide $ \O$ into $n$ intervals 
%\[
%\O_1^{(n)} =  \left [0, \frac{1}{n} \right [ , \O_2^{(n)} =  \left [\frac{1}{n}, \frac{2}{n} \right [,\ldots,\O_j^{(n)} =  \left [\frac{j-1}{n}, \frac{j}{n} \right [,\ldots,\O_n^{(n)} =  \left [ \frac{n-1}{n} , 1\right [,
%\]
%and let  $\mathcal{Q}_n$ denote the partition of $\O$, $\mathcal{Q}_n = \{ \O_i^{(n)}, i \in [n]\}$. Denote $\O_{ij}^{(n)} \eqdef \O_{i}^{(n)} \times \O_{j}^{(n)}$.

%\begin{itemize}[label = $\blacktriangleright$]
%\item

%\[
%g_{ni} = n  \int_{\O^{(n)}_{i} }  g(x)dx.
%\]
%\item 

% where $ (\Kn)_{ij}$ are obtained by averaging $K$ over the sets in $\mathcal{Q}_n$ 
%\begin{equation}
%\label{eq:kweight}
%(\Kn)_{ij} = n^2 \int_{\O^{(n)}_{i} \times \O^{(n)}_{j}}  K(x,y) dx dy.
%\end{equation}
%\end{itemize}

Our aim is to  study the relationship between the minimizer $\uast$ of $\Elam (\cdot, g, K)$ and the discrete minimizer $\unast$ of $\Elamn (\cdot, g_n, K_n)$ and estimate the error between solutions of discrete approximations and the solution of the continuum model. But the solution of problem~\eqref{def : maindiscpb} being discrete, it is convenient to introduce an intermediate model which is the continuum extension of the discrete solution. Towards this goal, we consider the piecewise constant injector $\injn$ of the discrete functions $\unast$ and $g_n$ into $L^2(\O)$, and of $\Kn$ into $L^{\infty}(\O^2)$, respectively. This injector $\injn$ is defined as
\begin{equation}\label{conitnuousextension}
\begin{aligned}
\injn u_n (x) &\eqdef \sum_{i=1}^{n} u_{ni} \chi_{\O^{(n)}_i}(x),\\
\injn g_n (x) &\eqdef \sum_{i=1}^{n} g_{ni} \chi_{\O^{(n)}_i}(x),\\
\injn \Kn (x,y)&\eqdef \sum_{i=1}^{n}  \sum_{j=1}^{n}K_{nij} \chi_{\O^{(n)}_i \times \O^{(n)}_j}(x, y) ,
\end{aligned}
\end{equation}
where we recall that $\chi_{\mathcal{C}}$ is the characteristic function of the set $\mathcal{C}$, i.e., takes $0$ on $\mathcal{C}$ and $1$ otherwise.

With these definitions, we have the following well-known properties whose proofs are immediate. We define the $\norm{\cdot}_{q,n}$ norm, for a given vector $u = (u_1,\cdots, u_n )^{\top} \in \R^n$, $q \in [1, +\infty[$, 
%we define the 
 %and  $ \norm{\cdot}_{p}$ norms as follows 
\[
\norm{u}_{q,n}=  \left ( \frac{1}{n} \sum\limits_{i=1}^{n} \aabs{u_i}^q\right )^{\frac{1}{q}}
\]
with the usual adaptation for $q=+\infty$.

\begin{lem}
\label{lem:projinj}
For a function $v \in L^q(\O)$, $q \in [1, + \infty]$, we have 
\begin{equation}
\norm{\projn v}_{q, n} \leq \norm{v}_{L^q(\O)};
\label{propprojector}
\end{equation}
and for $v_n \in \R^n$
\begin{equation}
\norm{\injn v_n}_{L^q(\O)} = \norm{v_n}_{q, n}.
\label{propinjector}
\end{equation}
In turn
\begin{equation}
\norm{\injn\projn v}_{L^q(\O)} \leq \norm{v}_{L^q(\O)} .
\label{proppvn}
\end{equation}
\end{lem}

It is immediate to see that the composition of the operators $\injn$ and $\projn$ yields the operator $\proj_{V_n}=\injn\projn$ which is the orthogonal projector on the subspace $V_n \eqdef \Span \ens{\chi_{\O^{(n)}_i}}{i \in [n]}$ of $L^1(\O)$.

%%%%%%%%%%%%%%%%%%%%%%%%%%%%%%%%%%%%%%%%%%%%%%%
\subsection{Graph limit theory}
\label{subsec:graphtheory}
We now briefly review some definitions and results from the theory of graph limits that we will need later since it is the key of our study of the discrete counterpart of the problem \eqref{def : mainpb} on dense deterministic graphs. We follow considerably~\cite{graph, lovs}, in which much more details can be found.

An undirected graph $G = \left ( V(G) , E(G) \right )$, where $V(G)$ stands for the set of nodes and $E(G) \subset V(G) \times V(G)$ denotes the edges set, without loops and parallel  edges is called simple.

Let $G_{n} = \pa{V(G_{n}) , E(G_{n})}$, $n \in \N^*$, be a sequence of dense, finite, and simple graphs, i.e; $\aabs{E(G_{n})} = O (\aabs{ V(G_{n})}^{2})$, where $\aabs{.}$ now denotes the cardinality of a set.

For two simple graphs $F$ and $G$, $\text{hom}(F , G)$ indicates the number of homomorphisms (adjacency-preserving maps) from $V (F )$ to $V (G)$. Then, it is worthwhile to normalize the homomorphism numbers and consider the homomorphism densities 
\[
t(F,G) = \frac{\text{hom} (F,G)}{\aabs{V(G)}^{\aabs{V(F)}}}.
\]
(Thus $t(F,G)$ is the probability that a random map of $V(F)$ into $V(G)$ is a homomorphism).
\begin{defi}\newcor{(cf.\cite{lovs})}
The sequence of graphs $\acc{G_{n}}_{n \in \N^*}$ is called convergent if $ t(F,G_{n})$ is convergent for every simple graph $F$.
\label{convergence}
\end{defi}
 
%\begin{rem}
%
%Note that $t(F,G_{n} ) = O(1)$ if $\aabs{E(G_{n})} = O (\aabs{ V(G_{n})}^{2})$  so that this definition is meaningful only for sequences of dense graphs. In the theory of graph limits, convergence in Definition~\ref{convergence} is called left-convergence. Since this is the only convergence of graph sequences that we use, we would refer to the left-convergent sequence as convergent (see \cite[Section~2.5]{borg}).
%\end{rem}

Convergent graph sequences have a limit object, which can be represented as a measurable \newcor{symmetric} function $K: \O^{2} \to \O$, here $\O$ stands for $[0,1]$. Such functions are called \textit{graphons}.

Let $\mathcal{K}$ denote the space of all bounded measurable functions  $K : \O^2 \to  \R$ such that $K(x,y)=K(y,x) $ for all $x,y \in [0,1]$. We also define $\mathcal{K}_{0} =\{ K \in \mathcal{K} : 0 \le K \le 1\}$ the set of all graphons.
\begin{prop}[{\cite[Theorem~2.1]{graph}}]
 For every convergent sequence of simple graphs, there is $K \in  \mathcal{K}_0$ such that 
 \begin{equation}
 t(F,G_{n}) \to t(F,K) \eqdef\int_{\O^{\aabs{V(F)}}}  \prod_{(i,j) \in E(F)} K(x_{i},x_{j}) dx
 \label{lim}
 \end{equation}
for every simple graph $F$. Moreover, for every $K \in \mathcal{K}_{0} $, there is a sequence of graphs $\acc{G_{n}}_{n \in \N^*}$ satisfying~\eqref{lim}.
\end{prop}
 
Graphon $K$ in \eqref{lim} which is uniquely determined up to measure-preserving transformations, is the limit of the convergent sequence $\acc{G_{n}}_{n \in \N^*}$.
%in the following sense: for every other limit function $K' \in \mathcal{K}_{0} $, there are measure-preserving map $\phi , \psi: \O \to \O$ such that $K(\phi(x), \phi(y)) = K' (\psi(x), \psi(y))$ (see \cite[Theorem~2.1]{graph}).
Indeed, every finite simple graph $G_{n}$ such that $V(G_{n}) = [n]$ can be represented by a function $K_{G_{n}} \in \mathcal{K}_{0}$ 
\begin{equation*}
K_{G_{n}} (x,y) = 
\begin{cases}
1 \quad \text{if}\quad (i,j) \in E(G_{n}) \quad \text{and} \quad (x,y) \in \O_{ij}^{(n)}, \\
0 \quad \text{otherwise}.
 \end{cases}
\end{equation*}
Hence, geometrically, the graphon $K$ can be interpreted as the limit of $K_{G_{n}}$ for the standard distance (called the cut-distance), see \cite[Theorem~2.3]{graph}.  An interesting consequence of this interpretation is that the space of graphs $G_{n}$, or equivalently pixel kernels $K_{G_{n}}$, is not closed under the cut distance. The space of graphons (larger than the space of graphs) defines the completion of this space.
%\[
% \norm{K}_{ \square} \eqdef\sup_{S,T \in \mathcal{L}_{\O}} \aabs{\int_{S \times T} K(x,y) dx dy},
%\]
% where $K  \in L^{1}( \O^2) $ and $\mathcal{L}_{\O}$ stands for the set of all Lebesgue measurable subsets of $\O$.
% Since for any $K  \in L^{1}( \O^2) $ 
%\[
%  \norm{K}_{ \square}   \le  \norm{K}_{ L^1(\O^2)},
%\]
%  convergence of $\{ K_{G_{n}}\}$ in the $L^1$-norm implies the convergence of the graph sequence $\{ G_{n}\}_n$ (\cite[Theorem~2.3]{graph}).

%%%%%%%%%%%%%%%%%%%%%%%%%%%%%%%%%%%%%%%%%%%%%%%%
\subsection{Lipschitz spaces}
\label{subsec: Lipspaces}
We introduce the Lipschitz spaces $\Lip(s,L^q(\O^d))$, for $d \in \{1,2\}$, $q \in [1, + \infty]$, which contain functions with, roughly speaking, $s$ "derivatives" in $L^q(\O^d)$~\cite[Ch.~2, Section~9]{devorelorentz93}. 
\begin{defi}\label{def:lipspaces}
For $F \in L^q(\O^d)$, $q \in [1,+\infty]$, we  define the (first-order) $L^q(\O^d)$ modulus of smoothness by
\begin{equation}
\omega(F,h)_q \eqdef
 \sup_{\bs z \in \R^d, |\bs z| < h} \pa{\int_{\bs x,\bs x +\bs z \in \O^d}\aabs{F(\bs x + \bs z)-F(\bs x)}^q d\bs x}^{1/q} .
\label{modsmooth}
\end{equation}
The Lipschitz spaces $\Lip(s,L^q(\O^d))$ consist of all functions $F$ for which
\[
\aabs{F}_{\Lip(s,L^q(\O^d))} \eqdef\sup_{h > 0} h^{-s} \omega(F,h)_q < +\infty .
\]
\end{defi}
We restrict ourselves to values $s \in ]0,1]$ since for $s > 1$, only constant functions are in $\Lip(s,L^q(\O^d))$. It is easy to see that $\aabs{F}_{\Lip(s,L^q(\O^d))}$ is a semi-norm. $\Lip(s,L^q(\O^d))$ is endowed with the norm
\[
\norm{F}_{\Lip(s,L^q(\O^2))} \eqdef \norm{F}_{L^q(\O^2)} +  \aabs{F}_{\Lip(s,L^q(\O^d))} .
\]
The space $\Lip(s,L^q(\O^2))$ is the Besov space $\mathbf{B}^s_{q,\infty}$~\cite[Ch.~2, Section~10]{devorelorentz93} which are very popular in approximation theory. In particular, $\Lip(s,L^{1/s}(\O^d))$ contains the space $\BV(\O^d)$ of functions of bounded variation on $\O^d$, i.e. the set of functions $F \in L^1(\O^d)$ such that their variation is finite:
\begin{equation*}
V_{\O^2}(F) \eqdef\sup_{h > 0}h^{-1}\sum_{i=1}^d\int_{\O^d}\aabs{F(\bs x + he_i)-F(\bs x)}d\bs x < + \infty, 
\end{equation*}
where $e_i, i \in \{1,d\}$ are the coordinate vectors in $\R^d$; see~\cite[Ch.~2, Lemma~9.2]{devorelorentz93}. Thus Lipschitz spaces are rich enough to contain functions with both discontinuities and fractal structure.

Let us define the piecewise constant approximation of a function $F \in L^q(\O^2)$ (a similar reasoning holds of course on $\O$) on a partition of $\O^2$ into cells 
\[
\O_{nij} \eqdef \ens{ ]x_{i-1}, x_i] \times ]y_{j-1}, y_j]} {(i,j) \in [n]^2}
\] 
of maximal mesh size $\deltan \eqdef \max\limits_{(i,j) \in [n]^2} \max(\aabs{x_{i} - x_{i-1}},\aabs{y_{j} - y_{j-1}})$,
\[
F_n(x,y) \eqdef \sum\limits_{i,j=1}^{n} F_{nij} \chi_{\O_{nij}}(x,y) , \quad F_{nij} = \frac{1}{\aabs{\O_{nij}}}\int_{\O_{nij}} F(x,y) dxdy .
\]
One may have recognized in these expressions non-equispaced versions of the projector and injector defined above.

We now state the following error bounds whose proofs use standard arguments from approximation theory; see~\cite[Section~6.2.1]{nonlocal} for details.
\begin{lem} 
There exists a positive constant $C_s$, depending only on $s$, such that for all $F \in \Lip(s,L^q(\O^d))$, $d \in \{1,2\}$, $s \in ]0,1]$, $q \in [1,+\infty]$,
\begin{equation}
\anorm{F - F_n}_{L^q(\O^d)} \le C_s \deltan^s \aabs{F}_{\Lip(s,L^q(\O^d))}.
\label{eq:lipspaceapprox}
\end{equation}
Let $p \in ]1, + \infty[$. If, in addition, $F \in L^{\infty}(\O^d)$, then \hl{
\begin{equation}
\anorm{F - F_n}_{L^p(\O^d)} \leq C\bpa{p,q,s,\anorm{F}_{L^{\infty}(\O^d)}}\deltan^{s \min\pa{1, q/p}} .
\label{eq: Lipg}
\end{equation}
where $C\bpa{p,q,s,\anorm{F}_{L^{\infty}(\O^d)}}$ is a positive constant depending only on $p$, $q$, $s$ and $ \anorm{F}_{L^{\infty}(\O^d)}$}.
\label{lem:spaceapprox}
\end{lem}

%\section{Existence and uniqueness of a minimizer}
\section{Well posedness}
\label{sec:exisvariational}

We start by proving existence and uniqueness of the minimizer for~\eqref{def : mainpb} and~\eqref{def : maindiscpb}. 
\begin{theo}
\label{theo:wellposed}
Suppose that $p \in [1, + \infty[$, $K$ is a nonnegative measurable mapping, and $g \in L^2(\O)$. Then, $\Elam(\cdot, g, K)$ has a unique minimizer in 
\[
\ens{u \in L^2(\O)}{\reg(u,K) \leq (2\lambda)^{-1}\norm{g}_{L^2(\O)}^2}, 
\]
and $\Elamn(\cdot, g_n, K_n)$ has a unique minimizer.
\end{theo}
\bpf{}
The arguments are standard (coercivity, lower semicontinuity and strict convexity) but we provide a self-contained proof (only for $\Elam(\cdot, g, K)$). \hl{First observe that from \cite[Proposition~9.32 and Proposition~9.5]{BauschkeCombBook11}, we infer that $\reg(\cdot,K)$ in \eqref{def : flow} is proper convex and lower semicontinuous}. Let $\acc{\uast_k}_{k \in \N}$ be a minimizing sequence in $L^2(\O)$. By optimality and Jensen's inequality, we have
\begin{equation}
\begin{aligned}
\norm{\uast_k}_{L^2(\O)}^2 
\leq 2\pa{2\lambda\Elam(\uast_k, g, K) + \norm{g}_{L^2(\O)}^2} 
&\leq 2\pa{2\lambda\Elam(0, g, K) + \norm{g}_{L^2(\O)}^2} \\
&= 4\norm{g}_{L^2(\O)}^2 < +\infty .
\end{aligned}
\label{eq:bnduastL2}
\end{equation}
Moreover
\begin{equation}
\reg(\uast_k,K) \leq \Elam(\uast_k, g, K) \leq \Elam(0, g, K) = \frac{1}{2\lambda}\norm{g}_{L^2(\O)}^2 < +\infty .
\label{eq:bnduastR}
\end{equation}
Thus $\norm{\uast_k}_{L^2(\O)}$ is bounded uniformly in $k$ so that the Banach-Alaoglu theorem for $L^2(\O)$ and compactness provide a weakly convergent subsequence (not relabelled) with a limit $\bar{u} \in L^2(\O)$. By lower semicontinuity of the $L^2(\O)$ norm and that of $\reg(\cdot,K)$, $\bar{u}$ must be a minimizer. The uniqueness follows from strict convexity of $\anorm{\cdot}^2_{L^2(\O)}$ and convexity of $\reg(\cdot,K)$.
\epf

\begin{rem}
Theorem~\ref{theo:wellposed} can be extended to linear inverse problems where the data fidelity in~$\Elam(0, g, K)$ is replaced by $\norm{g-\mathrm{A}u}^2_{L^2(\Sigma)}$, and where $\mathrm{A}$ is a continuous linear operator. The case where $\mathrm{A}: L^2(\O) \to L^2(\Sigma)$ is injective is immediate. The general case is more intricate and would necessitate appropriate assumptions on $\mathrm{A}$ and a Poincar\'e-type inequality. For instance, if $\mathrm{A}: L^p(\O) \to L^2(\Sigma)$, and the kernel of $\mathrm{A}$ intersects constant functions trivially, then using the Poincar\'e inequality in \cite[Proposition~6.19]{vaillo}, one can show existence and uniqueness in $L^p(\O)$, and thus in $L^2(\O)$ if $p \geq 2$. We omit the details here as this is beyond the scope of the paper.  
\end{rem}

We now turn to provide useful characterization of the minimizers $\uast$ and $\unast$. We stress that the minimization problem~\eqref{def : mainpb} that we deal with is considered over $L^2(\O)$ ($L^2(\O) \subset L^p(\O)$ only for $p \in [1,2]$) over which the function $\reg(\cdot,K)$ may not be finite (see \eqref{def : flow}). In correspondence, we will consider the subdifferential of the proper lower semicontinuous convex function $\reg(\cdot,K)$ on $L^2(\O)$ defined as
\[
\partial \reg(u,K) \eqdef \ens{\eta \in L^2(\O)}{\reg(v,K) \geq \reg(u,K) + \pds{\eta}{v-u}_{L^2(\O)}, ~~ \forall v \in L^2(\O)},
\]
and $\partial \reg(u,K) = \emptyset$ if $\reg(u,K) = +\infty$.

\begin{lem} 
Suppose that the assumptions of Theorem~\ref{theo:wellposed} hold. Then $\uast$ is the unique solution to~\eqref{def : mainpb} if and only if
\begin{equation}
\uast = \prox_{\lambda \reg(\cdot,K)}(g) \eqdef \pa{\Id+\lambda \partial \reg(\cdot,K)}^{-1}(g) .
\label{resolvent}
\end{equation}
Moreover, the proximal mapping $\prox_{\lambda \reg(\cdot,K)}$ is non-expansive on $L^2(\O)$, i.e., for $g_1, g_2 \in L^2(\O)$, the corresponding minimizers $\uast_1,\uast_2 \in L^2(\O)$ obey
\begin{equation}
\norm{\uast_1-\uast_2}_{L^2(\O)} \leq\norm{g_1 - g_2}_{L^2(\O)}.
\label{expansiveness}
\end{equation}
\label{lem:resolvent}
\end{lem}
A similar claim is easily obtained for~\eqref{def : maindiscpb} as well.

\bpf{}
The proof is again classical. By the first order optimality condition and since the squared $L^2(\O)$-norm is Fr\'echet differentiable, $\uast$ is the unique solution to~\eqref{def : mainpb} if, and only if the following monotone inclusion holds
\[
0 \in \uast-g + \lambda \partial \reg(\uast,K) .
\]
The first claim then follows. Writing the subgradient inequality for $\uast_1$ and $\uast_2$ we have
\begin{align*}
\reg(\uast_2,K) \geq \reg(\uast_1,K) +  \frac{1}{\lambda} \pds{g_1-\uast_1}{\uast_2-\uast_1}_{L^2(\O)} \\
\reg(\uast_1,K) \geq \reg(\uast_2,K) +  \frac{1}{\lambda} \pds{g_2-\uast_2}{\uast_1-\uast_2}_{L^2(\O)} .
\end{align*}
Adding these two inequalities we get
\[
\norm{\uast_2-\uast_1}_{L^2(\O)}^2 \leq \pds{\uast_2-\uast_1}{g_2-g_1}_{L^2(\O)},
\]
and we conclude upon applying Cauchy-Schwartz inequality.
\epf

We now formally derive the directional derivative of $\reg(\cdot,K)$ when $p \in ]1,+\infty[$. For this the symmetry assumption on $K$ is needed as well. Let $h \in L^2(\O)$. Then the following derivative exists
\[
\frac{d}{dt}\reg(u+th,K)\vert_{t=0} = \frac{1}{2}\int_{\O^2} K(x,y) \aabs{u(y)-u(x)}^{p-2}(u(y)-u(x))(v(y)-v(x)) dxdy.
\]
Since $K$ is symmetric, we apply the integration by parts formula in~\cite[Lemma~A.1]{nonlocal} (or split the integral in two terms and apply a change of variable $(x,y) \mapsto (y,x)$), to conclude that
\[
\frac{d}{dt}\reg(u+th,K)\vert_{t=0} = -\int_{\O^2} K(x,y) \aabs{u(y)-u(x)}^{p-2}(u(y)-u(x))v(x) dxdy = \pds{\A}{v}_{L^2(\O)} ,
\]
where 
\[
\A = -\int_{\O^2} K(x,y) \aabs{u(y)-u(x)}^{p-2}(u(y)-u(x)) dy
\] 
is precisely the nonlocal $p$-Laplacian operator, see~\cite{vaillo,nonlocal}. This shows that under the above assumptions, $\reg(\cdot,K)$ is Fr\'echet differentiable (hence G\^ateaux differentiable) on $L^2(\O)$ with Fr\'echet gradient $\A$.

\section{Error estimate for the discrete variational problem}
\label{sec:errorestimatevariational}

%%%%%%%%%%%%%%%%%%%%%%%%%%%%%%%%%%%%%%%%
\subsection{Main result}
\label{subsec:mainresult}
Our goal is to bound the difference between the unique minimizer of the continuum functional $\Elam(\cdot, g, K)$ defined on $L^2(\O)$ and the continuum extension by $\injn$ of that of $\Elamn (\cdot, g_n, \Kn)$. We are now ready to state the main result of this section.

\begin{theo} 
Suppose that $g \in L^{2}(\O)$ and $K$ is a nonnegative measurable, symmetric and bounded mapping. Let $\uast$ and $\unast$ be the unique minimizers of \eqref{def : mainpb} and \eqref{def : maindiscpb}, respectively. Then, we have the following error bounds.
\begin{enumerate}[label=(\roman*)]
\item If $p \in [1,2]$, then
\begin{equation}
\begin{aligned}
\norm{\In-\uast}_{L^2(\O)}^2 &\leq C \left(\norm{g-\Ig}_{L^2(\O)}^2 + \norm{g-\Ig}_{L^2(\O)} + \lambda\norm{K-\IK}_{L^{\frac{2}{2-p}}(\O^2)}\right.  \\
&+ \left. \lambda\norm{\uast - \IPn}_{L^{\frac{2}{3-p}}(\O)} \right),
\end{aligned}
\label{mainestimate1}
\end{equation}
where $C$ is a positive constant independent of $n$ and $\lambda$.
\item If $\inf_{(x,y) \in \O^2} K(x,y) \geq \kappa > 0$, then for any $p \in [1,+\infty[$,
\begin{equation}
\begin{aligned}
\norm{\In-\uast}_{L^2(\O)}^2 &\leq C \left(\norm{g-\Ig}_{L^2(\O)}^2 + \norm{g-\Ig}_{L^2(\O)} + \norm{K-\IK}_{L^\infty(\O^2)}\right.  \\
&+ \left. \lambda^{1/p}\norm{\uast - \IPn}_{L^p(\O)} \right),
\end{aligned}
\label{mainestimate2}
\end{equation}
where $C$ is a positive constant independent of $n$.
\end{enumerate}
\label{mainresult}
\end{theo}

A few remarks are in order before proceeding to the proof.
\hl{
\begin{rem}
{~}\\\vspace*{-0.5cm}
\begin{enumerate}[label=(\roman*)]
\item Observe that $2/(3-p) \leq p$ for $p \in [1,2]$. Thus by standard embeddings of $L^q(\O)$ spaces for $\O$ bounded, we have for $p \in [1,2]$
\[
\hspace*{-1cm}
\norm{K-\IK}_{L^{\frac{2}{2-p}}(\O^2)} \leq \norm{K-\IK}_{L^\infty(\O^2)} \qandq \norm{\uast - \IPn}_{L^{\frac{2}{3-p}}(\O)} \leq \norm{\uast - \IPn}_{L^p(\O)} ,
\]
which means that our bound in \eqref{mainestimate1} not only does not require an extra-assumption on $K$ but is also sharper than \eqref{mainestimate2}. The assumption on $K$ in the second statement seems difficult to remove or weaken. Whether this is possible or not is an open question that we leave to a future work.
\item We have made the dependence of the bound on $\lambda$ explicit on purpose. To see our motivation, assume that $g=\uorig+\varepsilon$, where $\uorig \in L^2(\O)$ is some true function and $\varepsilon \in L^2(\O)$ is some noise. Assume that $\partial \reg(\uorig,K) \neq \emptyset$, and let $\eta \in \partial \reg(\uorig,K)$, which is known in the inverse problem literature as a dual multiplier or certificate \cite{vaiterbookchap15}. Then
\[
\norm{\In-\uorig}_{L^2(\O)} \leq \norm{\In-\uast}_{L^2(\O)} + \norm{\uast-\uorig}_{L^2(\O)} .
\]
From \cite[Proposition~3.41]{scherzer2009variational}, one can show that
\[
\norm{\uast-\uorig}_{L^2(\O)} \leq 2\pa{\anorm{\varepsilon}_{L^2(\O)} + \lambda\anorm{\eta}_{L^2(\O)}} .
\]
With the standard choice $\lambda \sim \anorm{\varepsilon}_{L^2(\O)}$ we have $\anorm{\uast-\uorig}_{L^2(\O)}=O(\anorm{\varepsilon}_{L^2(\O)})$, and thus $\anorm{\uast-\uorig}_{L^2(\O)} \to 0$ as $\anorm{\varepsilon}_{L^2(\O)} \to 0$. Combining this with Theorem~\ref{mainresult} and the fact that
\[
\norm{g-\Ig}_{L^2(\O)} \leq \norm{\uorig-\injn\projn \uorig}_{L^2(\O)} + 2\norm{\varepsilon}_{L^2(\O)} ,
\]
one obtains an error bound of $\norm{\In-\uorig}_{L^2(\O)}$ as function of $\anorm{\varepsilon}_{L^2(\O)}$ and the discretization error of $\uorig$ and $K$. This error bound is dominated by that of $\uorig$ and $K$ as $\anorm{\varepsilon}_{L^2(\O)} \to 0$ fast enough. Having said this, as our focus here is on numerical consistency, in the rest of the paper, the dependence of the error bound on $\lambda$ will be absorbed in the constants.
\end{enumerate}
\end{rem}}

\bpf{}
\begin{enumerate}[label=(\roman*)]
\item Since $\Elam(\cdot, g, K)$ is a strongly convex function, we have 
\begin{equation}
\begin{aligned}
&\frac{1}{2\lambda}\norm{\In-\uast}_{L^2(\O)}^2 \\
&\leq \Elam(\In, g, K) - \Elam(\uast, g, K)\\
&\leq \bpa{\Elam(\In, g, K) - \Elamn(\unast, g_n, K_n)}
- \bpa{\Elam(\uast, g, K) -  \Elamn(\unast, g_n, K_n)}.
\end{aligned}
\label{ineq1 : strong}
\end{equation}
%Inequality \eqref{ineq : strong} implies that $  E(u, g, K) \leq E(\In, g, K)$. 
%\begin{equation}
%\begin{aligned}
%\frac{1}{2\lambda}\norm{\In-u}_{L^2(\O)}^2 &\leq E(\In, g, K) - E_n(u_n, g_n, K_n) - \left(E(u, g, K) -  E_n(u_n, g_n, K_n) \right)\\
%\label{ineq2 : strong}
%\end{aligned}
%\end{equation}
A closer inspection of $\Elam$ and $\Elamn$ and equality~\eqref{propinjector} allows to assert that
\begin{equation}
\Elam (\In, \Ig, \IK) = \Elamn (\unast, g_n, K_n).
\label{egalite}
\end{equation}
Now, applying the Cauchy-Schwarz inequality and using \eqref{egalite}, we have
\begin{equation}
\hspace*{-1cm}
\begin{aligned}
\Elam(\In, g, K) &= \frac{1}{2\lambda} \norm{\In - g}_{L^2(\O)}^2 + \reg (\In, K)\\
&= \frac{1}{2\lambda} \norm{\In - \Ig}_{L^2(\O)}^2 + \frac{1}{\lambda} \pds{\In - \Ig}{\Ig - g}_{L^2(\O)}\\
&+ \frac{1}{2\lambda} \norm{\Ig- g}_{L^2(\O)}^2 + \reg (\In, K)\\
&\leq \frac{1}{2\lambda} \norm{\In - \Ig}_{L^2(\O)}^2 + \frac{1}{\lambda} \norm{\In - \Ig}_{L^2(\O)}\norm{\Ig - g}_{L^2(\O)} \\
&+ \frac{1}{2\lambda} \norm{\Ig- g}_{L^2(\O)}^2 + \reg (\In, K)\\
&\leq \Elamn(\unast, g_n, K_n) +  \frac{1}{2\lambda} \norm{\Ig- g}_{L^2(\O)}^2 \\
&+ \frac{1}{\lambda} \norm{\In - \Ig}_{L^2(\O)}\norm{\Ig - g}_{L^2(\O)} +\bpa{\reg (\In, K) - \reg(\In, \IK)}\\
&\leq  \Elamn(\unast, g_n, K_n) +  \frac{1}{2\lambda} \norm{\Ig- g}_{L^2(\O)}^2 \\
&+ \frac{1}{\lambda} \norm{\In - \Ig}_{L^2(\O)}\norm{\Ig - g}_{L^2(\O)} \\
&+\frac{1}{2p} \aabs{\int_{\O^2} \bpa{K(x,y)- \IK(x,y)} \abs{\In(y) - \In(x)}^p dx dy} .
\end{aligned}
\label{ineq3 : strong}
\end{equation}
\hl{As we suppose that $g \in L^2(\O)$ and since $\unast$ is the (unique) minimizer of $\Elamn(\cdot, g_n, K_n)$, it is immediate to see, using \eqref{egalite} and \eqref{proppvn}, that
\begin{equation*}
\begin{aligned}
\frac{1}{2\lambda} \norm{\In - \Ig}_{L^2(\O)}^2 &\leq  \Elamn (\unast, g_n, K_n) \\
&\leq \Elamn(0, g_n, K_n)  = \Elam(0, \Ig, \IK)\\
&= \frac{1}{2\lambda}\norm{\Ig}_{L^2(\O)}^2
= \frac{1}{2\lambda}\norm{\injn\projn g}_{L^2(\O)}^2
\leq  \frac{1}{2\lambda}\norm{g}_{L^2(\O)}^2 < +\infty ,
%\leq  \frac{1}{2\lambda}\norm{g}_{L^{\infty}(\O)}^2< +\infty.
\end{aligned}
\end{equation*}}
and thus 
\begin{equation}
\norm{\In - \Ig}_{L^2(\O)} \leq \norm{g}_{L^2(\O)} \eqdef C_1 .
\label{ineq4 : stron}
\end{equation}
Since $p \in [1,2]$, by H\"older and triangle inequalities, and \eqref{eq:bnduastL2} applied to $\In$, we have that
\begin{equation}
\begin{aligned}
&\aabs{\int_{\O^2} \bpa{K(x,y)- \IK(x,y)} \abs{\In(y) - \In(x)}^p dx dy} \\
&\leq \norm{K - \IK}_{L^{\frac{2}{2-p}}(\O^2)}\pa{\int_{\O^2} \abs{\In(y) - \In(x)}^2 dx dy}^{p/2} \\
&\leq 2^{p}\norm{\In}_{L^2(\O)}^p \norm{K - \IK}_{L^{\frac{2}{2-p}}(\O^2)} \\ 
&\leq 2^{2p}\norm{\injn\projn g}_{L^2(\O)}^p \norm{K - \IK}_{L^{\frac{2}{2-p}}(\O^2)} \\
&\leq 2^{2p}\norm{g}_{L^2(\O)}^p \norm{K - \IK}_{L^{\frac{2}{2-p}}(\O^2)} = C_2 \norm{K - \IK}_{L^{\frac{2}{2-p}}(\O^2)},
\end{aligned}
\label{ineq5 : ston}
\end{equation}
where $C_2 \eqdef 2^{2p}C_1^p$.

We now turn to bounding the second term on the right-hand side of \eqref{ineq1 : strong}. Using~\eqref{proppvn} and the fact that $\unast$ is the (unique) minimizer of~\eqref{def : maindiscpb}, we have 
\begin{equation}
\hspace*{-0.5cm}
\begin{aligned}
   \Elam(\In, \Ig, \IK) &\leq   \Elam(\IPn, \Ig, \IK)\\
   &= \frac{1}{2 \lambda}\norm{\IPn - \injn \projn g}_{L^2(\O)}^2 + \reg(\IPn, \IK)\\
   &\leq  \frac{1}{2 \lambda}\norm{\uast- g}_{L^2(\O)}^2 +\reg(\uast, K)+ \reg(\IPn, \IK) - \reg(\uast, K)\\
   &\leq \Elam(\uast, g, K) + \left(\reg(\IPn, K) -\reg(\uast, K) \right)\\
   &+ (\reg(\IPn, \IK) - \reg(\IPn, K)).
\end{aligned}
\label{ineq4 : strong}
\end{equation}
We bound the second term on the right-hand side of \eqref{ineq4 : strong} by applying the mean value theorem on $[a(x,y),b(x,y)]$ to the function $t \in \R^+ \mapsto t^p$ with $a(x,y) = |\uast(y) - \uast(x)|$ and  $b(x,y) = |\IPn (y) - \IPn(x)|$.  Let $\eta(x,y) \eqdef \rho a(x,y) + (1-\rho) b(x,y)$, $\rho \in [0,1]$, be an intermediate value between $a(x,y)$ and $b(x,y)$. We then get
\begin{equation}
\begin{aligned}
&\abs{\reg(\IPn, K) -\reg(\uast, K)} \\
&= \abs{ \int_{\O^2} K (x,y) \pa{\abs{\IPn (y) - \IPn(x)}^p - \abs{\uast(y) - \uast(x)}^p} dx dy}\\
&= p\abs{ \int_{\O^2} K (x,y) \eta(x,y)^{p-1}\pa{\abs{\IPn (y) - \IPn(x)} - \abs{\uast(y) - \uast(x)}} dx dy}\\
&\leq  p C_3 \int_{\O^2} \eta(x,y)^{p-1}  \abs{\pa{\IPn (y) - \uast(y)} - \pa{\IPn (x) - \uast(x)}} dx dy \\
&\leq 2p C_3 \int_{\O^2} \eta(x,y)^{p-1}  \abs{\IPn (x) - \uast(x)} dx dy ,
\end{aligned}
\label{ineq1 : MVT}
\end{equation}
where we used the triangle inequality, symmetry after the change of variable $(x,y) \mapsto (y,x)$, and boundedness of $K$, say $\norm{K}_{L^\infty(\O^2)} \eqdef C_3$. Thus using H\"older and Jensen inequalities as well as~\eqref{proppvn}, and arguing as in \eqref{ineq5 : ston}, leads to
\begin{equation}
\begin{aligned}
&\abs{\reg(\IPn, K) -\reg(\uast, K)} \\
&\leq  2p C_3 \norm{\eta}_{L^2(\O^2)}^{p-1} \norm{\uast - \IPn}_{L^{\frac{2}{3-p}}(\O)} \\
&\leq  2p C_3 \pa{\rho\norm{a}_{L^2(\O^2)}+(1-\rho)\norm{b}_{L^2(\O^2)}}^{p-1} \norm{\uast - \IPn}_{L^{\frac{2}{3-p}}(\O)} \\
&\leq  2p C_3 \norm{a}_{L^2(\O^2)}^{p-1} \norm{\uast - \IPn}_{L^{\frac{2}{3-p}}(\O)} \\
&\leq  2^{2p-1}p C_3 \norm{g}_{L^2(\O)}^{p-1} \norm{\uast - \IPn}_{L^{\frac{2}{3-p}}(\O)} = C_4 \norm{\uast - \IPn}_{L^{\frac{2}{3-p}}(\O)} \\
\end{aligned}
\label{ineq6 : strong}
\end{equation}
where $C_4 \eqdef 2^{2p-1}pC_1^{p-1} C_3$.

To bound the last term on the right-hand side of \eqref{ineq4 : strong}, we follow the same steps as for establishing \eqref{ineq5 : ston} and get
\begin{equation}
\begin{aligned}
&|\reg(\IPn, \IK) - \reg(\IPn, K)| \\
&\leq \int_{\O^2} \abs{K(x,y)- \IK(x,y)} \abs{\IPn(y) - \IPn(x)}^p dx dy \\
&\leq C_2 \norm{K - \IK}_{L^{\frac{2}{2-p}}(\O^2)} .
\end{aligned}
\label{ineq7 : ston}
\end{equation}

Finally, plugging \eqref{ineq3 : strong}, \eqref{ineq4 : stron}, \eqref{ineq5 : ston}, \eqref{ineq4 : strong}, \eqref{ineq6 : strong} and \eqref{ineq7 : ston} into \eqref{ineq1 : strong}, we get the desired result.

\item The case $p \geq 2$ follows the same proof steps, except that now, we need to modify inequalities~\eqref{ineq5 : ston}, \eqref{ineq6 : strong} and \eqref{ineq7 : ston} which do not hold anymore. 

Under our assumption on $K$, and using~\eqref{eq:bnduastR}, \eqref{ineq5 : ston} now reads
\begin{equation}
\begin{aligned}
&\int_{\O^2} \abs{K(x,y) - \IK(x,y)} \abs{\In(y) - \In(x)}^p dx dy \\
&\leq \kappa^{-1} \norm{K - \IK}_{L^\infty(\O^2)} \int_{\O^2} \IK(x,y) \abs{\In(y) - \In(x)}^p dx dy \\
&= \kappa^{-1} \norm{K - \IK}_{L^\infty(\O^2)} \reg(\In,\IK) \\
&\leq (2\lambda\kappa)^{-1}C_1^2 \norm{K - \IK}_{L^\infty(\O^2)} ,
\end{aligned}
\label{ineq5ii : ston}
\end{equation}
where $C_1 = \norm{g}_{L^2(\O)}$ as in the proof of (i).

\hl{We embark from the last line of \eqref{ineq1 : MVT} to which we apply H\"older inequality and then Jensen inequality as well as~\eqref{proppvn} to get 
\begin{equation*}
\hspace*{-1.2cm}
\begin{aligned}
&\abs{\reg(\IPn, K) -\reg(\uast, K)} \\
&\leq  2p C_3 \norm{\eta}_{L^p(\O^2)}^{(p-1)} \norm{\uast - \IPn}_{L^p(\O)} \\
&\leq  2p C_3 \norm{a}_{L^p(\O^2)}^{(p-1)} \norm{\uast - \IPn}_{L^p(\O)} \\
&=  2p C_3 \pa{\int_{\O^2} \abs{\uast (y) - \uast(x)}^{p} dxdy}^{(p-1)/p} \norm{\uast - \IPn}_{L^p(\O)} .\end{aligned}
\end{equation*}
Now, by assumption on $K$ and using again~\eqref{eq:bnduastR}, we obtain
\begin{equation}
\hspace*{-1.2cm}
\begin{aligned}
&\abs{\reg(\IPn, K) -\reg(\uast, K)} \\
&\leq  2\kappa^{(1-p)/p}p C_3 \pa{\int_{\O^2} K(x,y) \abs{\uast (y) - \uast(x)}^{p} dxdy}^{(p-1)/p} \norm{\uast - \IPn}_{L^p(\O)} \\
&= 2\kappa^{(1-p)/p}p C_3 \pa{\reg(\uast,K)}^{(p-1)/p} \norm{\uast - \IPn}_{L^p(\O)} \\
&\leq 2(2\lambda\kappa)^{(1-p)/p}p C_3 C_1^{2(p-1)/p} \norm{\uast - \IPn}_{L^p(\O)} .
\end{aligned}
\label{ineq6ii : strong}
\end{equation}}

To get the new form of \eqref{ineq7 : ston}, we use~\eqref{proppvn}, \eqref{eq:bnduastR} and the assumption on $K$ to arrive at
\begin{equation}
\begin{aligned}
&|\reg(\IPn, \IK) - \reg(\IPn, K)| \\
&\leq \int_{\O^2} \abs{K(x,y)- \IK(x,y)} \abs{\IPn(y) - \IPn(x)}^p dx dy \\
&\leq \norm{K - \IK}_{L^\infty(\O^2)} \int_{\O^2} \abs{\uast(y) - \uast(x)}^p dx dy \\
&\leq \kappa^{-1}\norm{K - \IK}_{L^\infty(\O^2)} \int_{\O^2} K(x,y)\abs{\uast(y) - \uast(x)}^p dx dy \\
&= \kappa^{-1}\norm{K - \IK}_{L^\infty(\O^2)} \reg(\uast,K) \\
&\leq (2\lambda\kappa)^{-1}C_1^2\norm{K - \IK}_{L^\infty(\O^2)} .
\end{aligned}
\label{ineq7ii : ston}
\end{equation}
Plugging now \eqref{ineq3 : strong}, \eqref{ineq4 : stron}, \eqref{ineq4 : strong}, \eqref{ineq5ii : ston}, \eqref{ineq6ii : strong} and \eqref{ineq7ii : ston} into \eqref{ineq1 : strong}, we conclude the proof.
\end{enumerate}
\epf

%%%%%%%%%%%%%%%%%%%%%%%%%%%%%%%%%%%%%%%%
\subsection{Regularity of the minimizer}
The error bound of Theorem~\ref{mainresult} contain three terms: one which corresponds to the error in discretizing $g$, the second is the discretization error of the kernel $K$, and the last term reflects the discretization error of the minimizer $\uast$ of the continuum problem~\eqref{def : mainpb}. Thus, this form is not convenient to transfer our bounds to networks on graph and establish convergence rates. Clearly, we need a control on the term $ \norm{\IPn - \uast}_{L^q(\O)}$ on the right-hand side of~\eqref{mainestimate1}-\eqref{mainestimate2}. This is what we are about to do in the following key regularity lemma. In a nutshell, it states that if the kernel $K$ only depends on $|x-y|$ (as is the case for many kernels used in data processing), then as soon as the initial data $g$ belongs to some Lipschitz space, so does the minimizer $\uast$.

\begin{lem}
Suppose $g \in L^\infty(\O) \cap \Lip(s, L^q(\O))$ with $s\in ]0, 1]$ and $q \in [1, + \infty]$. Suppose furthermore that $K(x,y) = J(|x-y|)$, where $J$ is a nonnegative bounded measurable mapping on $\O$.
\begin{enumerate}[label=(\roman*)]
\item If $q \in [1,2]$, then $\uast \in \Lip(sq/2, L^q(\O))$. 
\item If $q \in [2,+\infty]$, then \hl{$\uast \in \Lip(s,L^2(\O))$}.
\end{enumerate}
\label{lem:lipgu}
\end{lem}

The boundedness assumption on $g$ can be removed for $q \geq 2$.

\bpf{}
We denote the torus $\tor \eqdef \R / 2 \Z$. For any function $u \in L^2(\O)$, we denote by $\bar{u} \in L^2(\tor)$ its periodic extension 
such that 
\begin{equation}
\baru(x)=
\begin{cases}
u(x) \quad \text{if} \quad x\in [0,1],\\
u(2-x) \quad \text{if} \quad x\in ]1,2],
\end{cases}
\label{def:extension}
\end{equation}
In the rest of the proof, we use letters with bars to indicate functions defined on $\tor$. 

Let us define
\[
\Ela(\barv, \barg, \barJ) \eqdef \frac{1}{\lambda} \norm{\barv- \barg}_{L^2(\tor)}^2 + \regb(\barv, \barJ)
\]
where
\[
\regb(\barv, \barJ) \eqdef \frac{1}{2p}\int_{\tor^2} \barJ(|x-y|) \abs{\barv(y)-\barv(x)}^p dx dy .
\]
Consider the following minimization problem 
%for example barvh, in the proof. Moreover we remind the reader for any discrete
%function fh defined on 
%h, the extended function Exth fh is also defined on 2
%h.
\begin{equation}
\min_{\barv \in L^2(\tor)}  \Ela(\barv, \barg, \barJ),
\label{def : periodicpb}
\end{equation}
which also has a unique minimizer by arguments similar to those of~Theorem~\ref{theo:wellposed}.
%One can easily see using the definitions of $\baru$, $\barg$ and $\barJ$ that 
%\[\displaystyle{\int_{\tor}} \baru(x) dx = 2 \displaystyle{\int_{\O}} u(x) dx \quad \text{and} \quad\displaystyle{\int_{\tor^2}} \barJ(|x-y|) dx dy = 4 \displaystyle{\int_{\O^2}} K(x,y) dx dy.
%\]
Since $\uast$ is the unique minimizer of \eqref{def : mainpb}, we have, using~\eqref{def:extension},
\begin{equation}
\begin{aligned}
\Ela(\ubast, \barg, \barJ) &= \frac{2}{\lambda} \norm{\uast - g}_{L^2(\O)}^2 + 4 \reg(\uast, J)\\
&= 4 \Elam(\uast, g, J) \\
&< 4 \Elam(v, g, J), \forall v \neq \uast \\
&=\Ela(\barv, \barg, \barJ), \forall \barv \neq \ubast ,
\end{aligned}
\end{equation}
which shows that $\ubast$ is the unique minimizer of \eqref{def : periodicpb}. Then, we have via Lemma~\ref{lem:resolvent}
%On the other hand, the unique minimizer of \eqref{def : mainpb} can be seen as the solution of the following problem 
%\begin{equation}
%u =\mathcal{J}_{\A}^{\lambda} g \eqdef \pa{\Id+\lambda \A}^{-1}g,
%\end{equation}
%with
%\[
%\A(u(x)) =  - \int_{\O} \K(x,y) \abs{u(y) - u(x) }^{p-2} (u(y) - u(x) ) dy,
%\]
%and the resolvent $\mathcal{J}_{\A}^{\lambda} $ is a single-valued non-expansive operator on $L^p(\O)$ since $\A$ is $m$-accretive (see~\cite{vaillo, Kato67} for more details). This observation leads us to 
\begin{equation}
\ubast = \prox_{\lambda/2 \regb(\cdot,\barJ)}(\barg).
\label{equa1}
\end{equation}
We define the translation operator
\[
(T_hv)(x) = v(x+h), \forall h \in \R .
\]
Now, using our assumption on the kernel $K$, that is $K(x,y) = J(|x-y|)$ (then invariant by translation), and periodicity of the functions on $\tor$, we have 
\begin{equation*}
\begin{aligned}
\Ela(\barv, T_h\barg, \barJ)&= \frac{1}{\lambda} \norm{\barv-T_h \barg}_{L^2(\tor)}^2 + \regb(\barv, \barJ)\\
&=  \frac{1}{\lambda} \norm{T_h(T_{-h}\barv-\barg)}_{L^2(\tor)}^2 \\
&+ \int_{\tor^2} \barJ (|x-y|)\abs{\barv((y+h)-h) - \barv((x+h)-h)}^p dx dy\\
&=  \frac{1}{\lambda} \norm{T_{-h}\barv-\barg}_{L^2(\tor)}^2 +  \int_{\tor^2} \barJ(|x-y|)\abs{T_{-h}\barv(y) - T_{-h}\barv(x)}^p dx dy\\
%&=  \frac{1}{\lambda} \norm{T_{-h}\barv-\barg}_{L^2(\tor)}^2 +  \int_{\tor^2} \barJ(|x-y|)\abs{T_{-h}\barv(y) - T_{-h}\barv(x)}^p dx dy\\
&=    \Ela(T_{-h}\barv, \barg, \barJ) .
\end{aligned}
\end{equation*}
This implies that the unique minimizer $\vbast$ of $ \Ela(\cdot, T_h\barg, \barJ)$ given by (see Lemma~\ref{lem:resolvent})
\begin{equation}
\vbast = \prox_{\lambda/2 \regb(\cdot,\barJ)}(T_h\barg) ,
\label{equa2}
\end{equation}
is also the unique minimizer of $\Ela(T_{-h}\cdot, \barg, \barJ)$.
%Relying on~\eqref{equa2} and using the fact that $\ubast$ is the unique minimizer of \eqref{def : periodicpb}, we have the following relationship between $\vbast$ and $\ubast$
%that the unique minimizer of problem~\eqref{equa2} satisfies
But since $ \Ela(\cdot, \barg, \barJ)$ has a unique minimizer $\ubast$, we deduce from~\eqref{equa1} and \eqref{equa2} that
\begin{equation}
T_h \prox_{\lambda/2 \regb(\cdot,\barJ)}(\barg) = \prox_{\lambda/2 \regb(\cdot,\barJ)}(T_h \barg) .
\label{equa3}	
\end{equation}
That is, the proximal mapping of $\lambda/2 \regb(\cdot,\barJ)$ commutes with translation.

We now split the two cases of $q$.
\begin{enumerate}[label=(\roman*)]
\item For $q \in [1,2]$: combining~\eqref{equa1}, \eqref{equa3}, \eqref{expansiveness}, \cite[Lemma~C.1]{nonlocal} and that $L^2(\O) \subset L^q(\O)$, we have
\begin{equation}
\begin{aligned}
\norm{T_h\ubast- \ubast}_{L^q(\tor)} 
&= \norm{\prox_{\lambda/2 \regb(\cdot,\barJ)}(T_h \barg)-\prox_{\lambda/2 \regb(\cdot,\barJ)}(\barg)}_{L^q(\tor)}\\
&\leq 2^{1/q-1/2}\norm{\prox_{\lambda/2 \regb(\cdot,\barJ)}(T_h \barg)-\prox_{\lambda/2 \regb(\cdot,\barJ)}(\barg)}_{L^2(\tor)} \\
&\leq 2^{1/q-1/2}\norm{T_h \barg - \barg}_{L^2(\tor)} \\
&\leq 2^{1/q-1/2}\bpa{2\norm{g}_{L^\infty(\O)}}^{1-q/2}\norm{T_h \barg - \barg}_{L^q(\tor)}^{q/2} \\
&= 2^{1/2}\bpa{\norm{g}_{L^\infty(\O)}}^{1-q/2} = C_1 \norm{T_h \barg - \barg}_{L^q(\tor)}^{q/2} .
\end{aligned}
\label{eq:bndwnonexp}
\end{equation}
Let $\O_h \eqdef \ens{x \in \O}{x+h \in \O}$.
Recalling the modulus of smoothness in \eqref{modsmooth}, we have 
\begin{equation}
\begin{aligned}
w(\uast,t)_q 
\eqdef \sup\limits_{\aabs{h} < t} \norm{T_{h}\uast -\uast }_{L^q(\O_h)} 
&\leq C_2 \sup\limits_{\aabs{h} < t} \norm{T_{h}\ubast -\ubast }_{L^q(\tor)} \\
&\leq C_1C_2 \pa{\sup\limits_{\aabs{h} < t} \norm{T_{h}\barg -\barg}_{L^q(\tor)}}^{q/2}\\
&=C_1C_2 w(\barg,t)_q^{q/2} \\
&\leq C_1C_2(C_3 w(g,t)_q)^{q/2}.
\end{aligned}
\end{equation}
We get the last inequality by applying the Whitney extension theorem~\cite[Ch.~6, Theorem~4.1]{devorelorentz93}.
Invoking Definition~\ref{def:lipspaces}, there exists a constant $C > 0$ such that
\begin{equation}
\aabs{\uast}_{\Lip(sq/2,L^q(\O))} \eqdef \sup\limits_{t>0} t^{-sq/2}w(\uast,t)_q \leq C \pa{\sup\limits_{t>0} t^{-s}w(g,t)_q}^{q/2} \leq C\aabs{g}_{\Lip(s,L^q(\O))}^{q/2} ,
\end{equation}
whence the claim follows after observing that $\uast \in L^2(\O) \subset L^q(\O)$.

\item \hl{For $q \in [2,+\infty]$, we combine~\eqref{equa1}, \eqref{equa3}, \eqref{expansiveness}, and that now $L^q(\O) \subset L^2(\O)$, to get
\[
\norm{T_h\ubast- \ubast}_{L^2(\tor)} \leq \norm{T_h \barg - \barg}_{L^2(\tor)} \leq 2^{1/2-1/q}\norm{T_h \barg - \barg}_{L^q(\tor)} .
\]
The rest of the proof is similar to that of (i).}
\end{enumerate}
\epf

In view of the regularity Lemma~\ref{lem:lipgu} and Theorem~\ref{mainresult}, one can derive convergence rates but only for $p \in [1,2]$. Indeed, the approximation bounds of Lemma~\ref{lem:spaceapprox} cannot be applied to $\uast - \IPn$ for $p \geq 2$ since the bound in Theorem~\ref{mainresult}(ii) is in the $L^p(\O)$ norm while Lemma~\ref{lem:lipgu} proves that $\uast$ is only in $\Lip(sq/2,L^2(\O))$. In particular, one cannot invoke \eqref{eq: Lipg} since there is no guarantee that $\uast$ is bounded. This is the reason why in the rest of the paper, we will only focus on the case $p \in [1,2]$.

\section{Application to dense deterministic graph sequences}
\label{sec:deteministicapplivariational}

The graph models we will consider here were used first in~\cite{medv} and then~\cite{nonlocal} to study networks on graphs for the evolution Cauchy problem, governed by the $p$-Laplacian in~\cite{nonlocal}. Throughout the section, we suppose that $p \in [1,2]$.

%%%%%%%%%%%%%%%%%%%%%%%%%%%%%%%%%%%%%%%%%%%%%%%%%%%%
\subsection{Networks on simple graphs}
\label{sec:simplevariational}
We first consider the case of a sequence of simple graphs converging to a $\{ 0,1\}$ graphon.
Briefly speaking, we define a sequence of simple graphs $G_n = \left (V(G_n),E(G_n) \right )$ such that $V(G_n) = [n]$ and 
\[
E(G_n) = \ens{(i,j) \in [n]^2}{ \O_{ij}^{(n)} \cap \cl{\supp(K)} \ne \emptyset},
\]
where $\cl{\supp(K)}$ is the closure of the support of $K$
\begin{equation}
\supp(K) = \ens{(x,y) \in \O^2}{K(x,y) \ne 0}.
\label{support}
\end{equation}
As we have mentioned in Section~\ref{subsec:graphtheory}, the kernel $K$ represents the corresponding graph limit, that is the limit as $n \to \infty$ of the function $K_{G_n} : \O^2 \to \{0,1\}$ such that 
\[
 K_{G_n}(x,y)   = \left\{\begin{array}{lr}
        1, \quad \text{if} \quad (i,j) \in E(G_n) \quad \text{and} \quad (x,y) \in \O_{ij}^{(n)},
        \\
      0  \quad \text{otherwise}.
        \end{array}
        \right.
\] 
As $n \to \infty $, $\acc{K_{G_n}}_{n \in \N^*}$ converges to the $\{0,1\}$-valued mapping $K$ whose support is defined by~\eqref{support}.
With this construction, the discrete counterpart of~\eqref{def : mainpb} on the graph $G_n$ is then given by 
\begin{equation}\tag{\textrm{$\mathcal{VP}_{s, n}^{\lambda, p}$}}
\min_{u_n \in \R^n}\acc{ \Elamn(u_n, g_n, K_n) \eqdef  \frac{1}{2 \lambda n } \norm{u_n - g_n}_2^2 + \frac{1}{2pn^2} \sum\limits_{(i,j) \in E(G_n)} \abs{u_{nj} - u_{ni}}^p},
\label{def: discretepbsimple}
\end{equation}
where the initial data $g_n$ is given by~\eqref{eq:averageg}.

For this model, $\IK(x,y)$ is the piecewise constant function such that for $(x,y) \in  \O_{ij}^{(n)}$, $(i,j) \in [n]^2$
\begin{equation}
\IK(x,y) = 
\begin{cases}
\frac{1}{|\O_{ij}^{(n)}|} \displaystyle{\int}_{\O_{ij}^{(n)}} K(x, y) dx dy   \quad \text{if} \quad  \O_{ij}^{(n)} \cap \cl{\supp(K)} \ne \emptyset, \\
0 \qquad \text{otherwise}.
\end{cases}
\label{eq:aveKsimple}
\end{equation}
%As $G_n $ is a simple graph, $\IK (\cdot, \cdot)$ is also a $\{0,1\}$-valued mapping.
Relying on what we did in~\cite{nonlocal}, the rate of convergence of the solution of the discrete problem to the solution of the limiting problem depends on the regularity of the boundary $\bd{\cl{\supp(K)}}$ of the closure of the support. \newcor{Following~\cite{medv}}, we recall the upper box-counting (or Minkowski-Bouligand) dimension of $\bd{\cl{\supp(K)}}$ as a subset of $\mathbb{R}^2$: 
\begin{equation}
\rho \eqdef \dim_B(\bd{\cl{\supp(K)}}) = \limsup_{\delta \rightarrow 0} \frac{\log N_{\delta}(\bd{\cl{\supp(K)}})}{- \log \delta},
\label{def: fractaldimension}
\end{equation}
where $ N_{\delta}(\bd{\cl{\supp(K)}}) $ is the number of cells of a $( \delta \times \delta )$-mesh that intersect $ \bd{\cl{\supp(K)}} $ 
(see \cite{falconer}).

\begin{theo} 
Assume that $p \in [1,2]$, $g\in L^{2}(\O)$. Let $\uast$ and $\unast$ be the unique minimizers of \eqref{def : mainpb} and  \eqref{def: discretepbsimple}, respectively. Then, the following hold.
\begin{enumerate}[label=(\roman*)]
\item We have 
\[
\norm{\In-\uast}_{L^2(\O)} \underset{n \to + \infty }{\longrightarrow} 0 .
\]
\label{convergence0}
\item For $p \in [1,2[$: assume moreover $g \in L^\infty(\O) \cap \Lip(s, L^q(\O))$, with $s  \in ]0,1]$ and $q \in [2/(3-p),2]$, that $\rho \in [0, 2[$ and that $K(x,y)=J(|x-y|)$, $\forall (x,y) \in \O^2$, with $J$ a nonnegative bounded measurable mapping on $\O$. Then for any $\epsilon >0$ there exists $N(\epsilon) \in \N$ such that for any $n \ge N(\epsilon)$
\begin{equation*}
\norm{\In - \uast}^2_{L^2(\O)} \le C n^{-\min\pa{sq/2,(2-p)(1-\frac{\rho+\epsilon}{2})}} ,
\end{equation*}
where $C$ is a positive constant independent of $n$.
%\item \todo{Check if we keep this or not.} For $p=2$: under the same assumptions as (ii), we have
%\begin{equation*}
%\norm{\In - \uast}^2_{L^2(\O)} \le C n^{-\min\pa{sq/2,2}} ,
%\end{equation*}
%where $C$ is a positive constant independent of $n$.
\end{enumerate}
\label{cor:simplegraphs}
\end{theo}

\bpf{}
\begin{enumerate}[label=(\roman*)]
\item In view of~\eqref{eq:averageg}, by the Lebesgue differentiation theorem (see e.g. \cite[Theorem~3.4.4]{pardoux}), we have 
\[
\Ig(x) \underset{n \to \infty}{\longrightarrow} g(x), \quad \IPn(x) \underset{n \to \infty}{\longrightarrow} \uast(x) \qandq \IK(x,y) \underset{n \to \infty}{\longrightarrow} K(x,y) 
\]
almost everywhere on $\O$ and $\O^2$, respectively. Combining this with Fatou's lemma and~\eqref{proppvn}, we have
\begin{align*}
\norm{g}_{L^2(\O)}^2 = \int_{\O} \aabs{\lim_{n \to \infty} \Ig(x)}^2 dx 
&= \int_{\O} \liminf_{n \to \infty} |\Ig(x)|^2 dx \\
&\leq \liminf_{n \to \infty} \norm{\Ig}_{L^2(\O)}^2 \\
&\leq \limsup_{n \to \infty} \norm{\injn\projn g}_{L^2(\O)}^2
\leq \norm{g}_{L^2(\O)}^2 ,
\end{align*}
which entails that $\lim_{n \to \infty} \norm{\Ig}_{L^2(\O)} = \norm{g}_{L^2(\O)}$. Similarly, we have 
\[
\lim_{n \to \infty} \norm{\IPn}_{L^\frac{2}{3-p}(\O)} = \norm{\uast}_{L^\frac{2}{3-p}(\O)}. 
\]
Since $g \in L^2(\O)$, $\uast \in L^2(\O) \subset L^{\frac{2}{3-p}}(\O)$ (Theorem~\ref{theo:wellposed}), we are in position to apply the Riesz-Scheff\'e lemma~\cite[Lemma~2]{riesz} to deduce that 
\[
\norm{\Ig - g}_{L^2(\O)} \underset{n \to \infty}{\longrightarrow} 0 \qandq \norm{\IPn - \uast}_{L^\frac{2}{3-p}(\O)} \underset{n \to \infty}{\longrightarrow} 0  . 
\]
Observe that for simple graphs, $\IK$ is not an orthogonal projection of $K$ (see~\eqref{eq:aveKsimple}) and thus, the above argument proof used for $g$ and $\uast$ does not hold. We argue however using the fact that $K$ is bounded, $|\O| < \infty$, and that $\forall n$ and $(x,y) \in \O^2$, $|\IK(x,y)| \leq \norm{K}_{L^\infty(\O)}$. We can thus invoke the dominated convergence theorem to get that 
\[
\norm{\IK - K}_{L^\frac{2}{2-p}(\O^2)} \underset{n \to \infty}{\longrightarrow} 0 .
\] 
Passing to the limit in~\eqref{mainestimate1}, we get the claim.

\item In the following $C$ is any positive constant independent of $n$. 
Since $g \in L^\infty(\O) \cap \Lip(s, L^q(\O))$, $q \leq 2$, and we are dealing with a uniform partition of $\O$ ($|\O_{i}^{(n)}|=1/n$, $\forall i \in [n]$), we get using inequality~\eqref{eq: Lipg} that 
\begin{equation}
\norm{\Ig - g}_{L^2(\O)} \leq C n^{-s\min\pa{1,q/2}} = C n^{-sq/2} .
\label{eq:estimategL2}
\end{equation}
By Lemma~\ref{lem:lipgu}(i), we have $\uast \in \Lip(sq/2,L^q(\O))$, and it follows from~\eqref{eq:lipspaceapprox} and the fact that $q \geq 2/(3-p)$ that
\begin{equation}
\norm{\IPn - \uast}_{L^{\frac{2}{3-p}}(\O)} \leq \norm{\IPn - \uast}_{L^q(\O)} \leq C n^{-sq/2} .
\label{eq:estimateuLq}
\end{equation}
Combining~\eqref{eq:estimategL2} and~\eqref{eq:estimateuLq}, we get 
\begin{equation}
\begin{aligned}
\norm{\Ig - g}_{L^2(\O)}^2 +  \norm{\Ig - g}_{L^2(\O)} + \norm{\IPn - \uast}_{L^{\frac{2}{3-p}}(\O)} 
&\leq C \bpa{n^{-sq} + n^{-sq/2}} \\
&\leq Cn^{-sq/2} .
\end{aligned}
\label{bornesuptermes}
\end{equation}
It remains to bound $\norm{K - \IK}_{L^{\frac{2}{2-p}}(\O^2)}$. For that, consider the set of discrete cells $\O_{ij}^{(n)}$ overlying the boundary of the support of $K$ 
\[
S(n) = \ens{(i,j) \in [n]^2}{\O_{ij}^{(n)} \cap \bd{\cl{\supp(K)}} \ne \emptyset} \qandq C(n) = \abs{S(n)}.
\]
For any $\epsilon > 0$ and sufficiently large $n$, we have 
\[
C(n) \le n^{\rho + \epsilon}.
\]
It is easy to see that $K$ and $\IK$ coincide almost everywhere on cells $\O_{ij}^{(n)}$ such that $(i,j) \notin S(n)$. Thus, for any  $\epsilon > 0$ and all sufficiently large $n$, we have 
\begin{equation}
\norm{K - \IK}_{L^{\frac{2}{2-p}}(\O^2)}^{\frac{2}{2-p}} \le C(n)n^{-2}  \le n^{-2(1- \frac{\rho + \epsilon}{2})}.
\label{bornesup}
\end{equation}
Inserting \eqref{bornesuptermes} and~\eqref{bornesup} into~\eqref{mainestimate1}, the desired result follows.

% \item For $p=2$, let $\O_{S(n)}=\bigcup_{(i,j) \in S(n)} \O_{ij}^{(n)}$. We then have
% \begin{align*}
% \norm{K - \IK}_{L^{\infty}(\O^2)} 
% &\leq \norm{K - \IK}_{L^{\infty}(\O^2 \setminus \O_{S(n)})} + \norm{K - \IK}_{L^{\infty}(\O_{S(n)})} \\
% &= \norm{K - \IK }_{L^{\infty}(\O_{S(n)})} \\
% &\leq \max_{(i,j) \in S(n)} \sup_{(x,y) \in \O_{ij}^{(n)}} |K(x,y) - \IK(x,y)| \leq n^{-2} .
% \end{align*}
\end{enumerate}
\epf

%%%%%%%%%%%%%%%%%%%%%%%%%%%%%%%%%%%%%%%%%%%%%%%%%%%%
\subsection{Networks on weighted graphs}
\label{sec:weightedvariational}
We now turn to the more general class of deterministic weighted graph sequences. The kernel $K$ is used to assign weights to the edges of the graphs considered below, we allow only positive weights. These weights $K_{nij}$ are obtained by averaging $K$ over the cells in the partition $\mathcal{Q}_n$ following~\eqref{eq:averageg}, and $\IK$ is given by~\eqref{conitnuousextension}.

Proceeding similarly to the proof of statement~\ref{convergence0} of Theorem~\ref{cor:simplegraphs}, we conclude immediately that
\[
\norm{\In-u}_{L^2(\O)} \underset{n \to + \infty }{\longrightarrow} 0.
\]
We are rather interested now in quantifying the rate of convergence in~\eqref{mainestimate1}. To do so, we need to add some regularity assumptions on the kernel $K$. 

\begin{theo}
Let $p \in [1,2[$, and assume that $g \in L^\infty(\O) \cap \Lip(s, L^q(\O))$, with $s  \in ]0,1]$ and $q \in [2/(3-p),2]$. Suppose moreover that $K(x,y)=J(|x-y|)$, $\forall (x,y) \in \O^2$, with $J$ a nonnegative bounded measurable mapping on $\O$. Let $\uast$ and $\unast$ be the unique minimizers of \eqref{def : mainpb} and  \eqref{def : maindiscpb}, respectively. Then, the following error bounds hold.
\begin{enumerate}[label=(\roman*)]
\item If $p \in [1,2[$, $K \in \Lip(s',L^{q'}(\O^2))$ and $(s',q') \in ]0,1] \times [1,+\infty[$, then 
\begin{equation}
\norm{\In-\uast}_{L^2(\O)}^2 \leq Cn^{-\min\pa{sq/2,s',s'q'(1-p/2)}} .
\label{ineq:weightedrate1}
\end{equation}
where $C$ is a positive constant independent of $n$.

In particular, if $g \in L^{\infty}(\O) \cap \BV(\O)$ and $K \in L^{\infty}(\O^2) \cap \BV(\O^2)$, then 
\begin{equation}
\norm{\In-u}_{L^2(\O)}^2 = O\bpa{n^{p/2-1}}.
\label{eq:weightedBVestimate1}
\end{equation}
\item If $p \in [1,2]$, $K \in \Lip(s',L^{q'}(\O^2))$ and $(s',q') \in ]0,1] \times [2/(2-p),+\infty]$, then
\begin{equation}
\norm{\In-\uast}_{L^2(\O)}^2 \leq Cn^{-\min\pa{sq/2,s'}} .
\label{ineq:weightedrate2}
\end{equation}
where $C$ is a positive constant independent of $n$.

In particular, if $g \in L^{\infty}(\O) \cap \BV(\O)$ then 
\begin{equation}
\norm{\In-u}_{L^2(\O)}^2 = O\bpa{n^{-\min\pa{1/2,s'}}} .
\label{eq:weightedBVestimate2}
\end{equation}
\end{enumerate}
\label{cor:weightedgraphs}
\end{theo}
\bpf{}
In the following $C$ is any positive constant independent of $n$. Under the setting of the theorem, for all cases, \eqref{bornesuptermes} still holds. It remains to bound $\norm{K - \IK}_{L^{\frac{2}{2-p}}(\O^2)}$. This is achieved using~\eqref{eq: Lipg} for case (i) and \eqref{eq:lipspaceapprox} for case (ii), which yields
\begin{equation}
\begin{cases}
\norm{K - \IK}_{L^{\frac{2}{2-p}}(\O^2)} \leq C n^{-s'\min\pa{1,q'(1-p/2)}} & \text{for case (i)} ,\\
\norm{K - \IK}_{L^{\frac{2}{2-p}}(\O^2)} \leq \norm{K - \IK}_{L^{q'}(\O^2)} \leq C n^{-s'} & \text{for case (ii)} .
\end{cases}
\label{inq:estimateK}
\end{equation}
Plugging \eqref{bornesuptermes} and~\eqref{inq:estimateK} into~\eqref{mainestimate1}, the bounds~\eqref{ineq:weightedrate1} and~\eqref{ineq:weightedrate2} follow.

We know that $\BV(\O) \subset \Lip(1/2,L^2(\O))$. Thus setting $s=s'=1/2$ and $q=q'=2$ in~\eqref{ineq:weightedrate1}, and observing that $1-p/2 \in [0,1/2]$, the bound~\eqref{eq:weightedBVestimate1} follows. That of~\eqref{eq:weightedBVestimate2} is immediate.
\epf

When $p=1$ (i.e., nonlocal total variation), $g \in L^{\infty}(\O) \cap \Lip(s,L^2(\O))$ and $K$ is a sufficiently smooth function, one can infer from Theorem~\ref{cor:weightedgraphs} that the solution to the discrete problem~\eqref{def : maindiscpb} converges to that of the continuum problem~\eqref{def : mainpb} at the rate $O(n^{-s})$. Moreover, if $g \in L^{\infty}(\O) \cap \BV(\O)$, then the best convergence rate is $O(n^{-1/2})$ which is attained precisely for $p=1$.
%This is to be compared to the convergence rate $O(n^{-s/(s+1)})$ established in~\cite[Theorem~4.1 and 5.1]{Wang2011} for the discretization of the local 2D ROF model.

\section{Application to random inhomogeneous graph sequences}
\label{sec:randomapplivariational}
We now turn to applying our bounds of Theorem~\ref{cor:weightedgraphs} to networks on random inhomogeneous graphs. %Throughout the section, it will be assumed implicitly that $K(x,y)=J(|x-y|)$, $\forall (x,y) \in \O^2$, with $J$ a nonnegative bounded measurable mapping on $\O$ to be able to invoke the regularity result of Lemma~\ref{lem:lipgu}. 

We start with the description of the random graph model we will use. This random graph model is motivated by the construction of inhomogeneous random graphs in~\cite{Ballobas2006, Ballobas2008}. It is generated as follows.

\begin{defi}
Fix $n \in \N^*$ and let $K$ be a symmetric measurable function on $\O^2$. Generate the graph $G_n = \pa{V(G_n),E(G_n)} \eqdef G_{q_n}(n,K)$ as follows: 
\begin{enumerate}[label=\arabic*)]
\item Generate $n$ independent and identically distributed (i.i.d.) random variables $\bX \eqdef (\bX_1, \cdots, \bX_n)$ from the uniform distribution on $\O$. Let $\acc{\bX_{(i)}}_{i = 1}^{n}$ be the order statistics of the random vector $\bX$, i.e. $\bX_{(i)}$ is the $i$-th smallest value.
\item Conditionally on $\bX$, join each pair $(i, j) \in [n]^2$ of vertices independently, with probability $q_n \wedgX$, i.e. for every $(i,j) \in [n]^2$, $i \neq j$,
\begin{equation}
 \P\pa{(i,j) \in E(G_n) | \bX} = q_n \wedgX , %= \min\pa{\frac{q_n}{\abs{\OXij}}\int_{\OXij}K(x, y) dx dy, 1},
 \label{def:graphmodel}
\end{equation}
where 
\begin{equation}
\wedgX \eqdef \min \pa{\frac{1}{\abs{\OXij}}\int_{\OXij}K(x, y) dx dy, 1/q_n},
\label{def:wedgX}
\end{equation}
and 
\begin{equation}
\OXij \eqdef ]\bX_{(i-1)}, \bX_{(i)}]  \times ]\bX_{(j-1)}, \bX_{(j)}] 
\end{equation}
where $q_n$ is nonnegative and uniformly bounded in $n$. 
\end{enumerate}
A graph $G_{q_n}(n,K)$ generated according to this procedure is called a $K$-random inhomogeneous graph generated by a random sequence $\bX$.
\label{def : randomgraph}
\end{defi}

We denote by $\bx = (\bx_1, \cdots, \bx_n)$  the realization of $\bX$. To lighten the notation, we also denote 
\begin{equation}
\O^{\bX}_{ni} \eqdef ]\bX_{(i-1)}, \bX_{(i)}], \quad \Oxi \eqdef   ]\bx_{(i-1)}, \bx_{(i)}], \tandt \Oxij \eqdef  ]\bx_{(i-1)}, \bx_{(i)}] \times   ]\bx_{(j-1)}, \bx_{(j)}] \quad  i,j \in [n].
\label{intervalpartition}
\end{equation}
As the realization of the random vector $\bX$ is fixed, we define 
\begin{equation}
\wedg \eqdef \min\pa{\frac{1}{\abs{\Oxij}}\displaystyle{\int_{\Oxij}}K(x, y) dx dy, 1/q_n}, \quad \forall (i,j) \in [n]^2, \quad  i \neq j.
\label{weightedgraph}
\end{equation}
In the rest of the paper, the following random variables will be useful. Let \linebreak$\Lam=\acc{\lam}_{(i,j) \in [n]^2, i\neq j}$ be a collection of independent random variables such that $q_n \lam$ follows a Bernoulli distribution with parameter $q_n \wedg$. We consider the independent random variables $\ups$ such that the distribution of $q_n \ups$ conditionally on $\bX=\bx$ is that of $q_n \lam$. Thus $q_n \ups$ follows a Bernoulli distribution with parameter $\EE\bpa{q_n \wedgX}$, where $\EE(\cdot)$ is the expectation operator (here with respect to the distribution of $\bX$). 

\medskip

\noindent We put the following assumptions on the parameters of the graph sequence $\acc{G_{q_n}(n,K)}_{n \in \N^*}$.
\begin{assum} We suppose that $q_n$ and $K$ are such that the following hold: 
\begin{enumerate}[label=\rm ({A}.\arabic*)]
\item $G_{q_n}(n,K)$ converges almost surely and its limit is the graphon $K \in L^{\infty} (\O^2)$;\label{assum:A1}
\item $\sup\limits_{n \geq 1} q_n < + \infty$.\label{assum:A2}
\end{enumerate}
\label{assumption}
\end{assum}

Graph models that verify \ref{assum:A1}-\ref{assum:A2} are discussed in \cite[Proposition~2.1]{nonlocalrandom}. 
They encompass the dense random graph model (i.e., with $\Theta(n^2)$ edges) extensively studied in~\cite{Lovaszdense2006,Lovaszrandom2011}, for which $q_n \geq c > 0$. This graph model allows also to generate sparse (but not too sparse) graph models; see~\cite{Ballobas2008}. That is graphs with $o(n^2)$ but $\omega(n)$ edges, i.e., that the average degree tends to infinity with $n$. For example, one can take $q_n=\exp(-\log(n)^{1-\delta})=o(1)$, where $\delta \in ]0,1[$.

%We now state our assumptions needed for the main results.
%\begin{assum} We suppose that the following hold : 
%\begin{enumerate}[label=\rm ({A}.\arabic*)]
%\item $G_n $ converges to a graphon $K \in L^{\infty} (\O^2)$;\label{assum:A1}
%\item $0 < \inf\limits_{n \geq 1} q_n \leq q_n \leq \sup\limits_{n \geq 1} q_n < + \infty$.\label{assum:A2}
%\end{enumerate}
%\label{assumptionvariational}
%\end{assum}
%We denote by $\bx = (\bx_1, \cdots, \bx_n)$  the realization vector of $\bX$. To lighten the notations, we  also denote 
%\begin{equation}
%\O^{\bX}_{ni} := ]\bX_{i-1}, \bX_i], \quad \Oxi :=   ]\bx_{i-1}, \bx_i], \quad \text{and} \quad   \Oxij :=  ]\bx_{i-1}, \bx_i] \times   ]\bx_{j-1}, \bx_j] \quad  i,j \in [n].
%\label{intervalpartitionvariational}
%\end{equation}
%As the realization of the random vector $X$ is fixed, we define the following parameter 
%\begin{equation}
%\wedg = \min \left \{\frac{1}{\abs{\Oxij}}\displaystyle{\int_{\Oxij}}K(x, y) dx dy, 1/q_n  \right \}, \quad (i,j) \in [n]^2, \quad  i \neq j.
%\label{weightedgraphvariational}
%\end{equation}
%In the rest of the paper, the following random variables will be useful.  We consider the i.i.d. random variables $\ups$ such that $q_n \ups$ follows a Bernoulli distribution with parameter $\EE\left (q_n \wedgX \right) $. Let $\lam$, $(i,j) \in [n]^2, i\neq j$, be i.i.d. random variables such that $q_n \lam $ follows a Bernoulli distribution with parameter $q_n \wedg$. It is immediate to see that the distribution of $\lam$ is that of $\ups$ conditionally to $\bX$.

%%%%%%%%%%%%%%%%%%%%%%%%%%%%%%%%%%%%%%%%%%%%%%%%%%%
\subsection{Networks on graphs generated by deterministic nodes}
\label{subsec: deterministicvariational}
In order to make our reasoning simpler, it will be convenient to assume first that the sequence $\bX$ is deterministic. Capitalizing on this result, we will then deal with the totally random model (i.e.; generated by random nodes) in Section~\ref{sec:applirandomvariational} by a simple marginalization argument combined with additional assumptions to get the convergence and quantify the corresponding rate. 
As we have mentioned before, we shall denote $\bx = (\bx_1, \cdots, \bx_n)$ as we assume that the sequence of nodes is deterministic. Relying on this notation, we define the parameter $\delta(n)$ as the maximal size of the spacings of $\bx$, i.e.,
\begin{equation}
\delta(n) = \max\limits_{i \in [n]} \abs{\bx_{(i)} - \bx_{(i-1)}}.
\label{spacingdeterministicvariational}
\end{equation}
Next, we consider the discrete counterpart of~\eqref{def : mainpb} on the graph $G_n$
\begin{equation}\tag{\textrm{$\mathcal{VP}_{d, n}^{\lambda, p}$}}
\min_{u_n \in \R^n}\acc{ \Elamn (u_n, g_n, K_n ) \eqdef \frac{1}{2 \lambda n } \norm{u_n - g_n}_2^2 + \frac{1}{2pn^2} \sum\limits_{i, j = 1}^{n}  \lam \abs{u_{nj} - u_{ni}}^p} ,
\label{def:discretepbramdom}
\end{equation}
where 
\[
g_i = \frac{1}{\abs{\Oxi}}\int_{\Oxi} g(x) dx.
\]
%Recall that $\lam$  are the i.i.d. random variables such that $q_n \lam $ follows a Bernoulli distribution with parameter $q_n \wedg$ defined above. To compare the solution of the continuum problem~\eqref{def : mainpb} to the solutions of~\eqref{def:discretepbramdom}, we define, using the injector $\injn$ defined in~\eqref{conitnuousextension}, the following intermediate random variable 
%\begin{equation}
%\Lam(x,y) \eqdef  \injn \lam (x, y) = \sum_{i=1}^{n}  \sum_{j=1}^{n}  \lam \chi_{\Oxij}(x, y)
%\label{continuousextensionrandom}
%\end{equation}
%and 
%\[
%\Ig(x) = \sum\limits_{i = 1}^{n} g_i  \chi_{\Oxi} (x).
%\]
%Next, we define the following step function 
%\[
%\Wedg(x,y) \eqdef \injn  \wedg (x,y)= \sum_{i=1}^{n}  \sum_{j=1}^{n}  \wedg \chi_{\Oxij}(x, y).
%\]

\begin{theo}
Suppose that $p \in [1,2[$, $g \in L^{2}(\O)$ and $K$ is a nonnegative measurable, symmetric and bounded mapping. Let $\uast$ and $\unast$ be the unique minimizers of \eqref{def : mainpb} and \eqref{def:discretepbramdom}, respectively. Let $p'=\frac{2}{2-p}$. 
\begin{enumerate}[label=(\roman*)]
\item There exist positive constants $C$ and $C_1$ that do not depend on $n$, such that for any $\beta > 0$
\begin{equation}
\hspace*{-0.5cm}
\begin{aligned}
\norm{\In-\uast}_{L^2(\O)}^2 &\leq C \left(\pa{\beta \frac{\log(n)}{n} + \frac{1}{q_n^{(p'-1)}n^{p'/2}}}^{1/p'} + \norm{g-\Ig}_{L^2(\O)}^2 \right. \\
&\left. + \norm{g-\Ig}_{L^2(\O)} + \norm{K-\injn\Wedgx}_{L^{p'}(\O^2)} + \norm{\uast - \IPn}_{L^{\frac{2}{3-p}}(\O)} \right),
\end{aligned}
\label{mainestimaterand1}
\end{equation}
with probability at least $1 - 2n^{-C_1 q_n^{2p'-1}\beta}$.

\item Assume moreover that $g \in L^\infty(\O) \cap \Lip(s, L^q(\O))$, with $s  \in ]0,1]$ and $q \in [2/(3-p),2]$, that $K(x,y)=J(|x-y|)$, $\forall (x,y) \in \O^2$, with $J$ a nonnegative bounded measurable mapping on $\O$, and $K \in \Lip(s',L^{q'}(\O^2))$, $(s',q') \in ]0,1] \times [p',+\infty]$ and $q_n\norm{K}_{L^\infty(\O^2)} \leq 1$. Then there exist positive constants $C$ and $C_1$ that do not depend on $n$, such that for any $\beta > 0$
\begin{equation}
\begin{aligned}
\norm{\In-\uast}_{L^2(\O)}^2 &\leq C \left(\pa{\beta \frac{\log(n)}{n} + \frac{1}{q_n^{(p'-1)}n^{p'/2}}}^{1/p'} +  \deltan^{-\min\pa{sq/2,s'}} \right),
\end{aligned}
\label{mainestimaterand2}
\end{equation}
with probability at least $1 - 2n^{-C_1 q_n^{2p'-1}\beta}$.
\end{enumerate}
\label{theo:deterministicnodesvarational}
\end{theo}

Before delving into the proof, some remarks are in order.
\begin{rem}
$~$\\\vspace*{-0.5cm}
\begin{enumerate}[label=(\roman*)]
\item The first term in the bounds~\eqref{mainestimaterand1}-\eqref{mainestimaterand2} can be replaced by
\[
\beta^{1/p'} \pa{\frac{\log (n)}{n}}^{1/p'} + \frac{1}{q_n^{(1-1/p')}n^{1/2}} .
\]
\item The last term in the latter bound can be rewritten as
\begin{equation}
n^{-1/2}q_n^{-(1-1/p')} = 
\begin{cases}
(q_n n)^{-1/2}				& \text{if } p' = 2 , \\
q_n^{1/p'}(q_n^2 n)^{-1/2} 		& \text{if } p' > 2 .
\end{cases}
\label{eq:varbndtermqn}
\end{equation}
Thus, if $\inf_{n \geq 1} q_n > 0$, as is the case when the graph is dense, then the term~\eqref{eq:varbndtermqn} is in the order of $n^{-1/2}$ with probability at least $1-n^{-c\beta}$ for some $c > 0$. If $q_n$ is allowed to be $o(1)$, i.e., sparse graphs, then \eqref{eq:varbndtermqn} is $o(1)$ if either $q_n n \to +\infty$ for $p'=2$, or $q_n^2 n \to +\infty$ for $p' > 2$. The probability of success is at least $1-e^{-C_1\beta\log(n)^{1-\delta}}$ provided that $q_n=\log(n)^{-\delta/(2p'-1)}$, with $\delta \in [0,1[$. All these conditions on $q_n$ are fulfilled by the inhomogenous graph model discussed above.

\item In fact, if $\inf_{n \geq 1} q_n \geq c > 0$, then we have $\sum_{n \geq 1} n^{-C_1q_n^{2p-1}\beta} \leq \sum_{n \geq 1} n^{-C_1c^{2p-1}\beta} < +\infty$ provided that $\beta > (C_1c^{2p-1})^{-1}$. Thus, if this holds, invoking the (first) Borel-Cantelli lemma, it follows that the bounds of Theorem~\ref{theo:deterministicnodesvarational} hold almost surely. The same reasoning carries over for the bounds of Theorem~\ref{theo:randomnodesvarational}.

\end{enumerate}
\label{rem:deterministicnodesvarational}
\end{rem}

\bpf{}
In the following $C$ is any positive constant independent of $n$. 
\begin{enumerate}[label=(\roman*)]
\item We start by arguing as in the proof of Theorem~\ref{mainresult}. Similarly to~\eqref{ineq1 : strong}, we now have
\begin{multline}
\frac{1}{2\lambda}\norm{\In-\uast}_{L^2(\O)}^2 
\leq \bpa{\Elam(\In, g, K) - \Elamn(\unast, g_n, \Lam)} \\
- \bpa{\Elam(\uast, g, K) -  \Elamn(\unast, g_n, \Lam)}.
\label{ineq1 : strongrand}
\end{multline}
The first term can be bounded similarly to~\eqref{ineq3 : strong}-\eqref{ineq4 : stron} to get
\begin{equation}
\begin{aligned}
&\Elam(\In, g, K) -  \Elamn(\unast, g_n, \Lam) \\
&\leq  C\Bigg(\norm{\Ig - g}_{L^2(\O)}^2 + \norm{\Ig - g}_{L^2(\O)} \\
&\quad + \aabs{\int_{\O^2} \bpa{K(x,y)- \injn\Lam(x,y)} \abs{\In(y) - \In(x)}^p dx dy}\Bigg) \\
&\leq C\Bigg(\norm{\Ig - g}_{L^2(\O)}^2 + \norm{\Ig - g}_{L^2(\O)} \\
&\quad + \aabs{\int_{\O^2} \bpa{K(x,y)- \IWedg(x,y)} \abs{\In(y) - \In(x)}^p dx dy} \\
&\quad + \aabs{\int_{\O^2} \bpa{\IWedg(x,y)- \injn\Lam(x,y)} \abs{\In(y) - \In(x)}^p dx dy}\Bigg) .
\end{aligned}
\label{ineq3 : strongrand}
\end{equation}
The second term in~\eqref{ineq3 : strongrand} is $O\pa{\norm{K - \IWedg}_{L^{p'}(\O^2)}}$, see \eqref{ineq5 : ston}. For the last term, we have using Jensen and H\"older inequalities,
\begin{equation}
\begin{aligned}
&\aabs{\int_{\O^2} \bpa{\IWedg(x,y)- \injn\Lam(x,y)} \abs{\In(y) - \In(x)}^p dx dy} \\
&\leq 2^{p-1} \Bigg(\int_{\O}\aabs{\int_{\O} \bpa{\IWedg(x,y)- \injn\Lam(x,y)} dy} \abs{\In(x)}^p dx \\
&\quad + \int_{\O}\aabs{\int_{\O} \bpa{\IWedg(x,y)- \injn\Lam(x,y)} dx} \abs{\In(y)}^p dy\Bigg) \\
&\leq C \Bigg(\pa{\int_{\O}\aabs{\int_{\O} \bpa{\IWedg(x,y)- \injn\Lam(x,y)} dy}^{p'} dx}^{1/p'} \\
&\quad + \pa{\int_{\O}\aabs{\int_{\O} \bpa{\IWedg(x,y)- \injn\Lam(x,y)} dx}^{p'} dy}^{1/p'}\Bigg) \\
&= C \pa{\norm{Z_{n}}_{p',n} + \norm{W_{n}}_{p',n}} ,
\end{aligned}
\label{eq:bndlamKrand}
\end{equation}
where 
\[
Z_{ni} \eqdef \frac{1}{n} \sum_{j=1}^n \pa{\wedg-\lam} \qandq W_{nj} \eqdef \frac{1}{n} \sum_{i=1}^n \pa{\wedg-\lam} .
\]
By virtue of~\cite[Lemma~A.1]{nonlocalrandom}, together with \ref{assum:A2} and the fact that $p' \geq 2$, there exists a positive constant $C_1$, such that for any $\beta > 0$
\[
\P\pa{\anorm{Z_n}_{p',n} \ge \varepsilon } \le n^{-C_1 q_n^{2p'-1}\beta} ,
\]
%\[
%\P\pa{\anorm{Z_n}_{p',n} \ge \varepsilon } \le n^{-C_1 q_n^{2p'-1}\beta} = n^{-C_1 q_n^{2p'-1}} ,
%\]
with
\begin{equation}
\label{eq:epsilon}
\varepsilon =  \pa{\beta \frac{\log(n)}{n}+\frac{1}{q_n^{(p'-1)}n^{p'/2}}}^{1/p'} .
\end{equation}
The same bound also holds for $\norm{W_{n}}_{p',n}$. A union bound then leads to
\begin{equation}
\norm{Z_{n}}_{p',n} + \norm{W_{n}}_{p',n} \leq 2\varepsilon
\label{eq:bndZW}
\end{equation}
with probability at least $1-2n^{-C_1 q_n^{2p'-1}\beta}$.

Let us now turn to the second term in~\eqref{ineq1 : strongrand}. Using~\eqref{proppvn} and the fact that $\unast$ is the unique minimizer of~\eqref{def:discretepbramdom}, we have 
\begin{equation}
\begin{aligned}
   \Elam(\In, \Ig, \injn\Lam) - \Elam(\uast, g, K) 
   &\leq \pa{\reg(\IPn, K) -\reg(\uast, K)}\\
   &+ \bpa{\reg(\IPn, \IWedg) - \reg(\IPn, K)} \\
   &+ \bpa{\reg(\IPn, \injn\Lam) - \reg(\IPn, \IWedg)} .
\end{aligned}
\label{ineq4 : strongrand}
\end{equation}
The first term is bounded as in~\eqref{ineq6 : strong}, which yields
\begin{equation}
\begin{aligned}
\abs{\reg(\IPn, K) -\reg(\uast, K)} \leq  C \norm{\uast - \IPn}_{L^{\frac{2}{3-p}}(\O)}.
\end{aligned}
\label{ineq6 : strongrand}
\end{equation}
The second term follows from~\eqref{ineq7 : ston}
\begin{equation}
\begin{aligned}
\abs{\reg(\IPn, \IWedg) - \reg(\IPn, K)} \leq C \norm{K - \IWedg}_{L^{p'}(\O^2)} .
\end{aligned}
\label{ineq7 : stonrand}
\end{equation}
The last term is upper-bounded exactly as in~\eqref{eq:bndlamKrand} and \eqref{eq:bndZW}.

Inserting \eqref{ineq3 : strongrand}, \eqref{eq:bndlamKrand}, \eqref{eq:bndZW}, \eqref{ineq4 : strongrand}, \eqref{ineq6 : strongrand} and \eqref{ineq7 : stonrand} into \eqref{ineq1 : strongrand}, we get the claimed bound.

\item Insert \eqref{bornesuptermes} and~\eqref{inq:estimateK} into~\eqref{mainestimaterand1} after replacing $1/n$ by $\deltan$.
\end{enumerate}
\epf

%%%%%%%%%%%%%%%%%%%%%%%%%%%%%%%%%%%%%%%%%%%%%%%%
\subsection{Networks on graphs generated by random nodes}
\label{sec:applirandomvariational}
Let us turn now to the totally random model. The discrete counterpart of~\eqref{def : mainpb} on the totally random sequence of graphs $\acc{G_{q_n}}_{n \in \N^*}$ is given by 
\begin{equation}\tag{\textrm{$\mathcal{VP}_{r, n}^{\lambda, p}$}}
\min_{u_n \in \R^n}\acc{ \Elamn (u_n, g_n, K_n ) \eqdef  \frac{1}{2 \lambda n } \norm{u_n - g_n}_2^2 + \frac{1}{n^2} \sum\limits_{i, j = 1}^{n}  \ups\abs{u_{nj} - u_{ni}}^p},
\label{def:discretepbtotalrandom}
\end{equation}
where we recall that the random variables $\ups$ are independent with $q_n \ups$ following the Bernoulli distribution with parameter $\EE\left (q_n \wedgX \right) $ defined above. 

Observe that for the totally random model, $\deltan$ is a random variable. Thus, we have to derive a bound on it. In~\cite[Lemma~3.2]{nonlocalrandom}, it was shown that
\begin{equation}
\deltan \leq  t \frac{\log (n)}{n},
\label{eq:diconcentrationdeltan}
\end{equation}
with probability at least $1 - n^{-t}$, where $t \in ]0,e[$.

Combining this bound with Theorem~\ref{theo:deterministicnodesvarational} (after conditioning and integrating) applied to the totally random sequence $\acc{G_{q_n}}_{n \in \N^*}$, we get the following result.
\begin{theo}
Suppose that $p \in [1,2[$, $g \in L^{2}(\O)$ and $K$ is a nonnegative measurable, symmetric and bounded mapping. Let $\uast$ and $\unast$ be the unique minimizers of \eqref{def : mainpb} and \eqref{def:discretepbtotalrandom}, respectively. Let $p'=\frac{2}{2-p}$. 
\begin{enumerate}[label=(\roman*)]
\item There exist positive constants $C$ and $C_1$ that do not depend on $n$, such that for any $\beta > 0$
\begin{equation}
\hspace*{-0.5cm}
\begin{aligned}
\norm{\In-\uast}_{L^2(\O)}^2 &\leq C \left(\pa{\beta \frac{\log(n)}{n} + \frac{1}{q_n^{(p'-1)}n^{p'/2}}}^{1/p'} + \norm{g-\Ig}_{L^2(\O)}^2 \right. \\
&\left. + \norm{g-\Ig}_{L^2(\O)} + \norm{K-\injn\WedgX}_{L^{p'}(\O^2)} + \norm{\uast - \IPn}_{L^{\frac{2}{3-p}}(\O)} \right),
\end{aligned}
\label{mainestimaterandnodes1}
\end{equation}
with probability at least $1 - 2n^{-C_1 q_n^{2p'-1}\beta}$.

\item Assume moreover that $g \in L^\infty(\O) \cap \Lip(s, L^q(\O))$, with $s  \in ]0,1]$ and $q \in [2/(3-p),2]$, that $K(x,y)=J(|x-y|)$, $\forall (x,y) \in \O^2$, with $J$ a nonnegative bounded measurable mapping on $\O$, that $K \in \Lip(s',L^{q'}(\O^2))$, $(s',q') \in ]0,1] \times [p',+\infty]$ and $q_n\norm{K}_{L^\infty(\O^2)} \leq 1$. Then there exist positive constants $C$ and $C_1$ that do not depend on $n$, such that for any $\beta > 0$ and $t \in ]0,e[$
\begin{equation}
\begin{aligned}
\norm{\In-\uast}_{L^2(\O)}^2 &\leq C \left(\pa{\beta \frac{\log(n)}{n} + \frac{1}{q_n^{(p'-1)}n^{p'/2}}}^{1/p'} +  \pa{t\frac{\log(n)}{n}}^{\min\pa{sq/2,s'}} \right),
\end{aligned}
\label{mainestimaterandnotes2}
\end{equation}
with probability at least $1 - \bpa{2n^{-C_1 q_n^{2p'-1}\beta}+n^{-t}}$.
\end{enumerate}
\label{theo:randomnodesvarational}
\end{theo}

\bpf{}
Again, $C$ will be any positive constant independent of $n$.
\begin{enumerate}[label=(\roman*)]
\item Let
\begin{align*}
\varepsilon' 
&= C \left(\pa{\beta \frac{\log(n)}{n}+C \frac{1}{q_n^{(p'-1)}n^{p'/2}}}^{1/p'} + \norm{g-\Ig}_{L^2(\O)}^2 + \norm{g-\Ig}_{L^2(\O)} \right. \\
&\qquad \left. + \norm{K-\injn\WedgX}_{L^{p'}(\O^2)} + \norm{\uast - \IPn}_{L^{\frac{2}{3-p}}(\O)} \right) .
\end{align*}
Using~\eqref{mainestimaterand1}, and independence of this bound from~$\bx$, we have
\begin{equation*}
\begin{aligned}
\P\pa{\norm{\In-\uast}_{L^2(\O)}^2  \ge \varepsilon'} 
&= \frac{1}{\abs{\O}^n}\int_{\O^n} \P\pa{\norm{\In-\uast}_{L^2(\O)}^2  \ge \varepsilon'| \bX = \bx} d\bx \\
&\leq \frac{1}{\abs{\O}^n}\int_{\O^n} 2n^{-C_1 q_n^{2p'-1}\beta} d\bx\\
&=2n^{-C_1 q_n^{2p'-1}\beta}.
\end{aligned}
\end{equation*}

\item Recall $\varepsilon$ in~\eqref{eq:epsilon} and $\kappa = C\pa{t\frac{\log(n)}{n}}^{\min\pa{sq/2,s'}}$. 
Denote the event
{\small
\[
\hspace*{-0.75cm}
A_1: \acc{\norm{g-\Ig}_{L^2(\O)}^2 + \norm{g-\Ig}_{L^2(\O)} + \norm{K-\injn\WedgX}_{L^{p'}(\O^2)} + \norm{\uast - \IPn}_{L^{\frac{2}{3-p}}(\O)} \leq \kappa}.
\]}
In view of \eqref{bornesuptermes}, \eqref{inq:estimateK} and~\eqref{eq:diconcentrationdeltan}, and that under our assumptions $\WedgX=K_{n}^{\bX}$, we have
\[
\P\pa{A_1} \geq \P\pa{\deltan \leq t \frac{\log(n)}{n}} \geq 1-n^{-t} .
\]
Let the event 
\begin{align*}
A_2&: \acc{\norm{Z_{n}}_{p',n} + \norm{W_{n}}_{p',n} \leq 2\varepsilon} ,
\end{align*}
and denote $A_i^c$ the complement of the event $A_i$. It then follows from~\eqref{eq:bndZW} and the union bound that 
\begin{align*}
\P\pa{\norm{\In-\uast}_{L^2(\O)}^2 \leq 2C\varepsilon + \kappa} 
&\geq \P\pa{A_1\cap A_2} = 1 - \P\pa{ A_1^c \cup A_2^c} \\
&\geq 1 - \sum_{i=1}^2 \P\pa{A_i^c} \geq 1 - \pa{2n^{-C_1 q_n^{2p'-1}\beta} + n^{-t}} ,
\end{align*}
which leads to the claimed result.
\end{enumerate}
\epf

When $p=1$ (i.e., nonlocal total variation), $g \in L^{\infty}(\O) \cap \Lip(s,L^2(\O))$ and $K$ is a sufficiently smooth function, one can deduce from Theorem~\ref{theo:randomnodesvarational} that with high probability, the solution to the discrete problem~\eqref{def:discretepbtotalrandom} converges to that of the continuum problem~\eqref{def : mainpb} at the rate $O\pa{\pa{\frac{\log(n)}{n}}^{-\min(1/2,s)}}$. Compared to the deterministic graph model, there is overhead due to the randomness of the graph model which is captured in the rate and the extra-logarithmic factor.

\section{Numerical results}
\label{sec:numerical}

In this section, we will apply the variational regularization problem~\eqref{def : maindiscpb} to a few applications, and illustrate numerically our bounds.

\subsection{Minimization algorithm}
The algorithm we will describe in this subsection is valid for any $p \in [1,+\infty]$\footnote{Obviously $\lim_{p \to +\infty} \frac{1}{p} \norm{\cdot}_p^p = \iota_{\ens{u_n \in \R^n}{\anorm{u_n}_{\infty} \leq 1}}$.}. The minimization problem~\eqref{def : maindiscpb} can be rewritten in the following form
\begin{equation}
\min_{u_n \in \R^n} \frac{1}{2} \norm{u_n - g_n}_2^2 + \frac{\lambda_n}{p}\norm{\nabla_{K_n} u_n}^p_{p} ,
\label{def:maindiscpbprimal}
\end{equation}
where $\lambda_n=\lambda/(2n)$, $\nabla_{K_n}$ is the (nonlocal) weighted gradient operator with weights $K_{nij}$, defined as
\begin{align*}
\nabla_{K_n}:
& \R^n \to \R^{n \times n} \\
& u_n \mapsto V_n, ~~ V_{nij} = K_{nij}^{1/p}(u_{nj}-u_{ni}), \quad \forall (i,j) \in [n]^2 .
\end{align*}
This is a linear operator whose adjoint, the (nonlocal) weighted divergence operator denoted $\div_{K_n}$. It is easy to show that
\begin{align*}
\div_{K_n}:
& \R^{n \times n} \to \R^{n} \\
& V_n \mapsto u_n, ~~ u_{ni} = \sum_{m=1}^n K_{nmi}^{1/p} V_{nmi} - \sum_{j=1}^n K_{nij}^{1/p} V_{nij}, \quad \forall i \in [n] .
\end{align*}
Problem~\eqref{def:maindiscpbprimal} can be easily solved using standard duality-based first-order algorithms. For this we follow~\cite{FadiliTV10}. 

By standard conjugacy calculus, the Fenchel-Rockafellar dual problem of~\eqref{def:maindiscpbprimal} reads
\begin{equation}
\min_{V_n \in \R^{n \times n}} \frac{1}{2} \norm{g_n - \div_{K_n}V_n}_2^2 + \frac{\lambda_n}{q}\norm{V_n/\lambda_n}^q_{q} ,
\label{def:maindiscpbdual}
\end{equation}
where $q$ is the H\"older dual of $p$, i.e. $1/p+1/q=1$. One can show with standard arguments that the dual problem~\eqref{def:maindiscpbdual} has a convex compact set of minimizers for any $p \in [1,+\infty[$. Moreover, the unique solution $\unast$ to the primal problem~\eqref{def:maindiscpbprimal} can be recovered from any dual solution $\Vnast$ as
\[
\unast = g_n - \div_{K_n}\Vnast .
\]
It remains now to solve~\eqref{def:maindiscpbdual}. The latter can be solved with the (accelerated) FISTA iterative scheme~\cite{Nesterov83,fista2009,chambolle2015convergence} which reads in this case 
\begin{equation}
\begin{aligned}
W_n^{k} 	&= V_n^{k} + \frac{k-1}{k+b}(V_n^{k} - V_n^{k-1}) \\
V_n^{k+1}	&= \prox_{\gamma\frac{\lambda_n}{q}\anorm{\cdot/\lambda_n}^q_{q}}\pa{W_n^{k} + \gamma \nabla_{K_n}\bpa{g_n-\div_{K_n}(W_n^{k})}} \\
u_n^{k+1} 	&=  g_n - \div_{K_n}V_n^{k+1} ,
\end{aligned}
\label{eq:fistadual}
\end{equation}
where $\gamma \in \big]0,\bpa{\sup_{\anorm{u_n}_2 = 1}\anorm{\nabla_{K_n}u_n}_2}^{-1}\big]$, $b > 2$, and we recall that $\prox_{\tau F}$ is the proximal mapping of the proper lsc convex function $F$ with $\tau > 0$, i.e.,
\[
\prox_{\tau F}(W) = \uArgmin{V \in \R^{n \times n}} \frac{1}{2}\norm{V-W}_2^2 + \tau F(V) .
\] 
The convergence guarantees of scheme~\eqref{eq:fistadual} are summarized in the following proposition.
\begin{prop}
\label{prop:convfistadual}
The primal iterates $u_n^{k}$ converge to $\unast$, the unique minimizer of~\eqref{def : maindiscpb}, at the rate
\[
\norm{u_n^{k} - \unast}_2 = o(1/k) .
\]
\end{prop}

\bpf{}
Combine~\cite[Theorem~2]{FadiliTV10} and \cite[Theorem~1.1]{AttouchFISTA15}.
\epf

Let us turn to the computation of the proximal mapping $\prox_{\gamma\frac{\lambda_n}{q}\anorm{\cdot/\lambda_n}^q_{q}}$. Since $\anorm{\cdot}^q_{q}$ is separable, one has that
\[
\prox_{\gamma\frac{\lambda_n}{q}\anorm{\cdot/\lambda_n}^q_{q}}(W) = \pa{\prox_{\gamma\frac{\lambda_n}{q}\aabs{\cdot/\lambda_n}^q}(W_{ij})}_{(i,j) \in [n]^2} .
\] 
Moreover, as $\aabs{\cdot}^q$ is an even function on $\R$, $\prox_{\gamma\frac{\lambda_n}{q}\aabs{\cdot/\lambda_n}^q}$ is an odd mapping on $\R$, that is,
\[
\prox_{\gamma\frac{\lambda_n}{q}\aabs{\cdot/\lambda_n}^q}(W_{ij}) = \prox_{\gamma\frac{\lambda_n}{q}\aabs{\cdot/\lambda_n}^q}(|W_{ij}|)\sign\pa{W_{ij}} .
\]
In a nutshell, one has to compute $\prox_{\gamma\frac{\lambda_n}{q}\aabs{\cdot/\lambda_n}^q}(t)$ for $t \in \R^+$. We distinguish different situations depending on the value of $q$:
\begin{itemize}
\item $q=+\infty$ (i.e., $p=1$): this case amounts to computing the orthogonal projector on $[-\lambda_n,\lambda_n]$, which reads
\[
t \in \R^+ \mapsto \proj_{[-\lambda_n,\lambda_n]}(t) = \min\bpa{t,\lambda_n} .
\]

\item $q=1$ (i.e., $p=+\infty$): this case corresponds to the well-known soft-thresholding operator, which is given by
\[
t \in \R^+ \mapsto \prox_{\gamma |\cdot|}(t) = \max\bpa{t-\gamma,0}  .
\]

\item $q=2$ (i.e., $p=2$): it is immediate to see that
\[
\prox_{\gamma/(2\lambda_n) |\cdot|^2}(t) = \frac{t}{1+\gamma/\lambda_n} .
\]

\item $q \in ]1,+\infty[$: in this case, as $|\cdot|^q$ is differentiable, the proximal point $\prox_{\gamma\frac{\lambda_n}{q}\aabs{\cdot/\lambda_n}^q}(t)$ is the unique solution $\alpha^\star$ on $\R^+$ of the non-linear equation
\[
\alpha - t + \gamma\alpha^{p-1}/\lambda_n = 0 .
\]
\end{itemize}

\subsection{Experimental setup}
We apply the scheme~\eqref{eq:fistadual} to solve~\eqref{def:maindiscpbprimal} in two applicative settings with nonlocal regularization on (weighted) graphs. The first one pertains to denoising of a function defined on a 2D point cloud, and the second one to signal denoising. In the first setting, the nodes of the graph are the points in the cloud and $u_{ni}$ is the value of point/vertex index $i$. For signal denoising, each graph node correspond  to a signal sample, and $u_{ni}$ is the signal value at node/sample index $i$. We chose the nearest neighbour graph with the standard weighting kernel $e^{-|\bf{x}-\bf{y}|}$ when $|\bf{x}-\bf{y}| \leq \delta$ and $0$ otherwise, where $\bf{x}$ and $\bf{y}$ are the 2D spatial coordinates of the points for the point cloud\footnote{For the 2D case, $(\bf{x},\bf{y})$ are not to be confused with the "coordinates" $(x,y)$ of the graphon on the continuum, though there is a bijection from one to another.}, and sample index for the signal case. 

\begin{figure}[ht]
\centering
\includegraphics[width=0.5\textwidth]{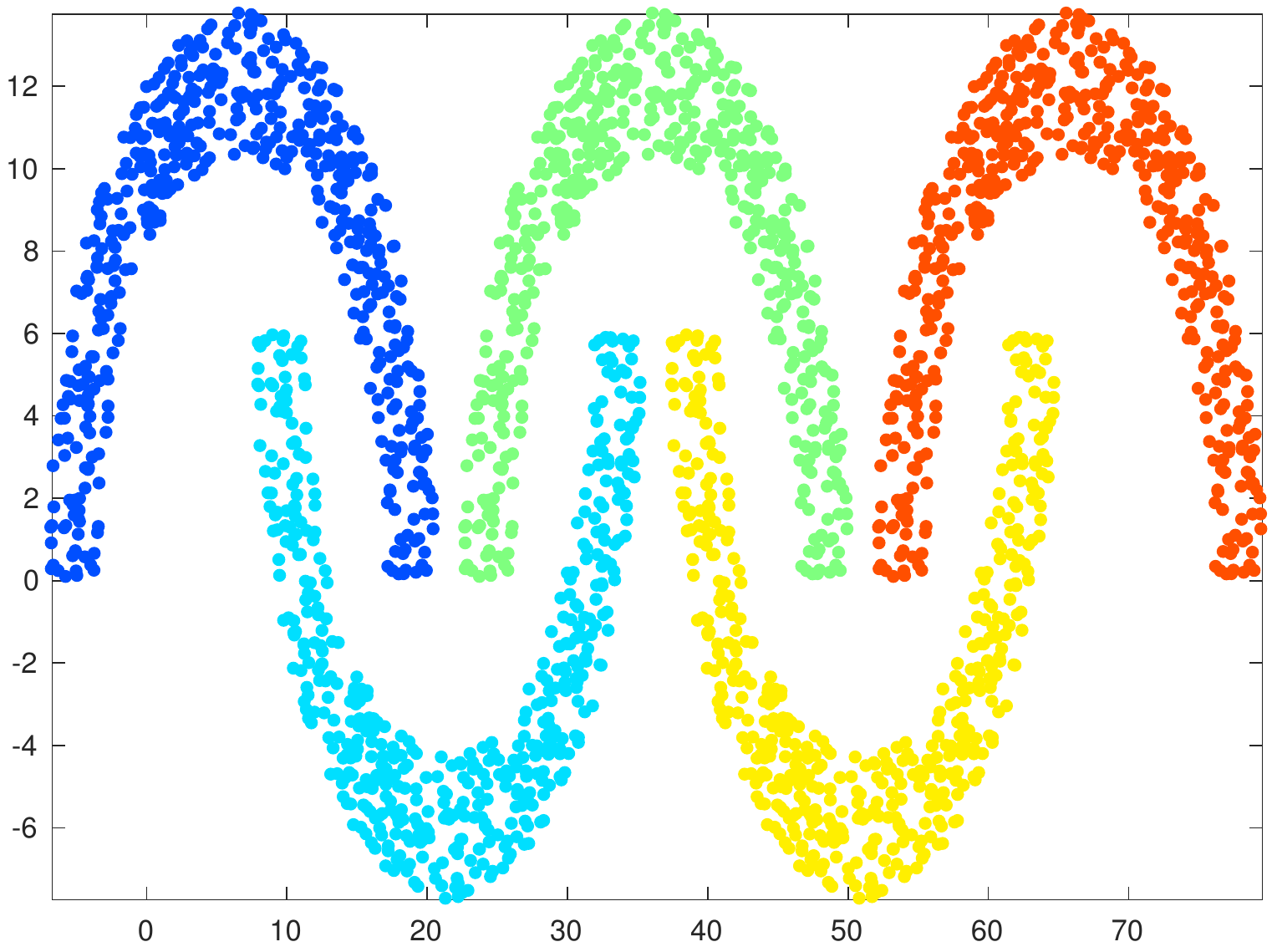}
\caption{Original point cloud with $N=2500$ points.}
\label{fig:pointcloudoriginal}
\end{figure}

\begin{figure}[ht]
\centering
\includegraphics[width=\textwidth]{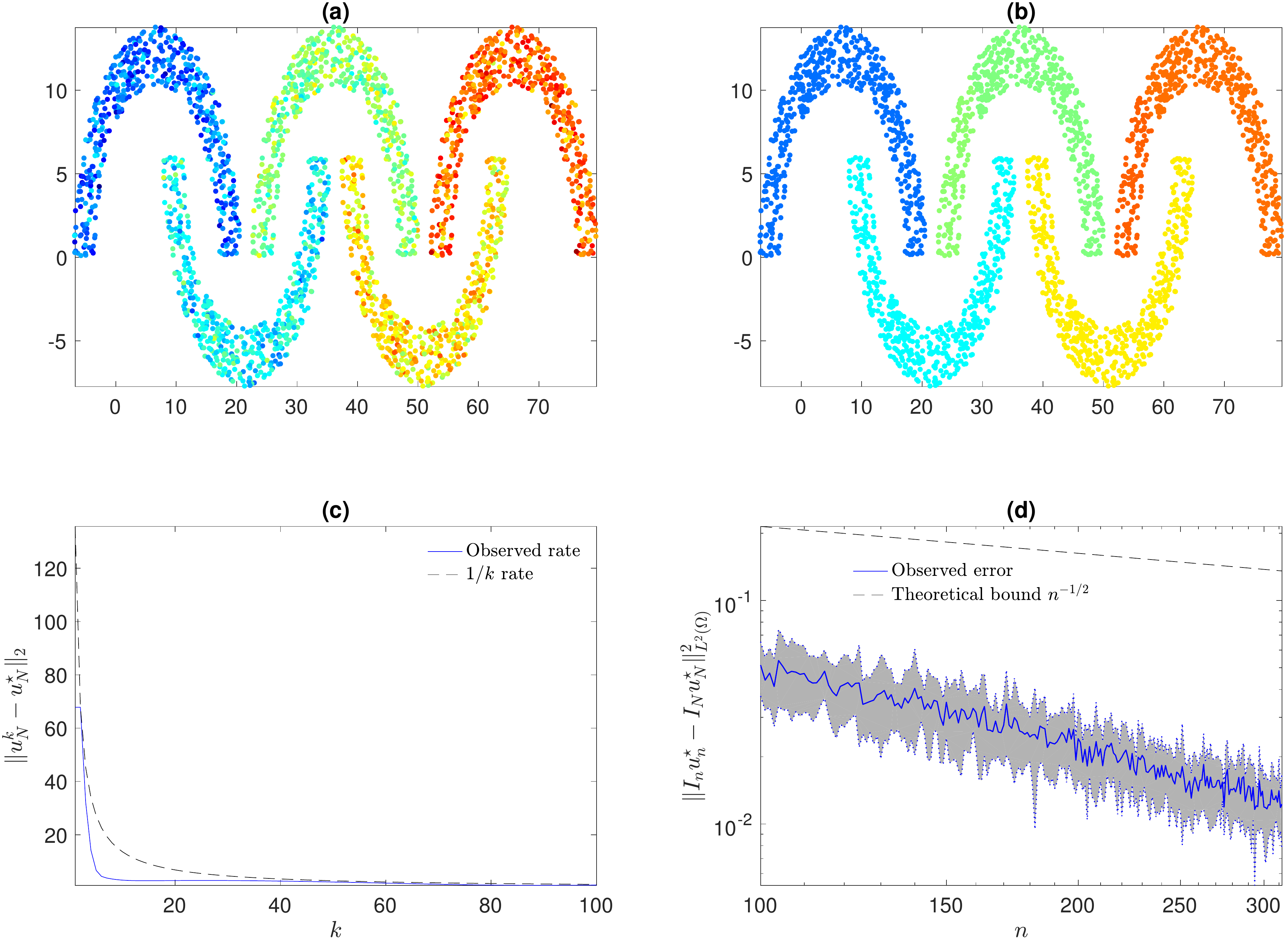} 
\caption{Results for point cloud denoising with $p=1$. (a) Noisy point cloud. (b) Recovered point cloud by solving~\eqref{def:maindiscpbprimal}. (c) Primal convergence criterion $\norm{u_n^k-\unast}_2$ as a function of the iteration counter $k$. (d) Mean error $\norm{\In - I_N\uast_N}^2_{L^2(\O)}$ across replications as a function of $n$.}
\label{fig:pointcloudresultp1}
\end{figure}

\paragraph{Application to point cloud denoising}
The original point cloud used in our numerical experiments is shown in~Figure~\ref{fig:pointcloudoriginal}. It consists of $N=2500$ points that are not on a regular grid. The function on this point cloud, denoted $u^0_{N}$, is piecewise-constant taking 5 values (5 clusters) in $[5]$. A noisy observation $g_N$ (see Figure~\ref{fig:pointcloudresultp1}(a)) is then generated by adding a white Gaussian noise of standard deviation $0.5$ to $u^0_{N}$. Given the piecewise-constancy of $u^0_{N}$, we solved~\eqref{def:maindiscpbprimal} with the natural choice $p=1$. The result is shown in Figure~\ref{fig:pointcloudresultp1}(b). Figure~\ref{fig:pointcloudresultp1}(c) displays the evolution of $\norm{u_N^k-u_N^\star}_2$ as a function of the iteration counter $k$, which confirms the theoretical rate $o(1/k)$ predicted above.

To illustrate our consistency results, $\uast$ is needed while it is not known in our case. Therefore, we argue as follows. We consider the continuum extension of $I_N u_N^\star$ as a reference and compute $\norm{I_n \unast-I_N u_N^\star}_{L^2(\O)}$ for varying $n \ll N$. \hl{By the triangle inequality, $\norm{I_n \unast-I_N u_N^\star}_{L^2(\O)}$ is clearly dominated by $\norm{I_n \unast-\uast}_{L^2(\O)}$}. Thus, for each value of $n \in [100,N/8]$, $n$ nodes are drawn uniformly at random in $[N]$ and $g_n$ is generated, which is a sampled version of $g_{N}$ at those nodes. This is replicated $20$ times. For each replication, we solve~\eqref{def:maindiscpbprimal} with $g_n$ and the same regularization parameter $\lambda$, and we compute the mean across the 20 replications of the squared-error $\norm{\In-I_N u_N^\star}_{L^2(\O)}^2$. The result is depicted in Figure~\ref{fig:pointcloudresultp1}(d). The gray-shaded area corresponds to one standard deviation of the error over the 20 replications. One indeed observes that the average error decreases at a rate consistent with the $O(n^{-1/2})$ predicted by our results (see discussion after Theorem~\ref{cor:weightedgraphs} with $s=1/2$).

\begin{figure}[ht]
\centering
\includegraphics[width=\textwidth]{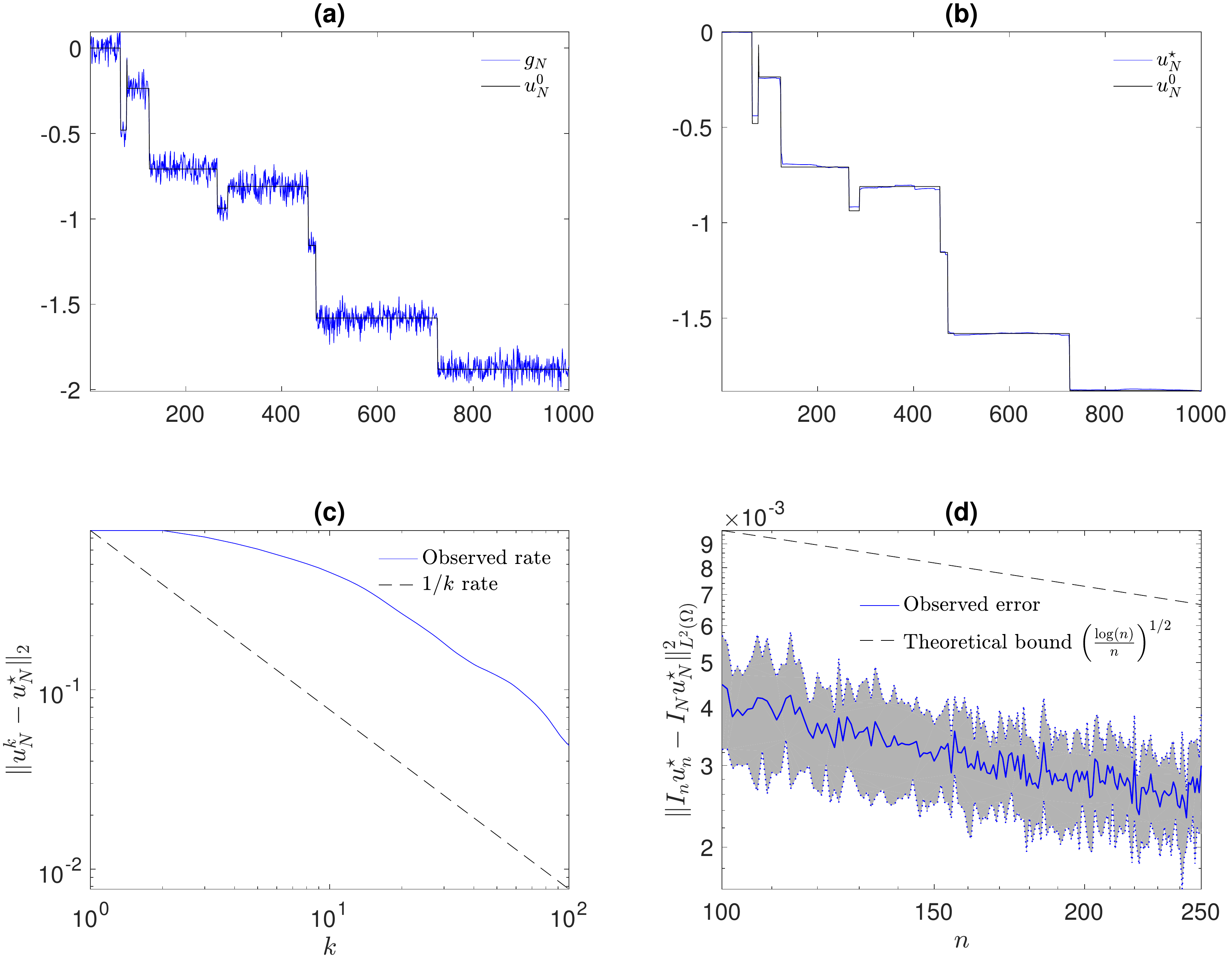} 
\caption{Results for signal denoising with $p=1$. (a) Noisy and original signal. (b) Denoised and original signal foor $N=1000$. (c) Primal convergence criterion $\norm{u_n^k-\unast}_2$ as a function of the iteration counter $k$. (d) Mean error $\norm{\In - I_N\uast_N}^2_{L^2(\O)}$ as a function of $n$.}
\label{fig:signalresultp1}
\end{figure}

\paragraph{Application to signal denoising}
In this experiment, we choose a piecewise-constant signal shown in Figure~\ref{fig:signalresultp1}(a) for $N=1000$ together with its noisy version $g_N$ with additive white Gaussian noise of standard deviation $0.05$. Figure~\ref{fig:signalresultp1}(b) depicts the denoised signal $u_N^\star$ by solving~\eqref{def:maindiscpbprimal} with $p=1$ and hand-tuned $\lambda$. Figure~\ref{fig:signalresultp1}(c) also confirms the $o(1/k)$ rate predicted above on $\norm{u_N^k-u_N^\star}_2$.

We now illustrate the consistency bound result on a random sequence of graphs \linebreak$\acc{G_{q_n}(n,K)}_{n \in [100,N/4]}$ generated according to Definition~\ref{def : randomgraph} with $q_n=1$. For each value of $n \in [100,N/4]$, $n$ nodes are drawn uniformly at random in $[N]$, and $g_n$ is generated, which is a sampled version of $g_{N}$ at those nodes. $n^2$ independent Bernoulli variables $\lam$ each with parameter $K_{nij}$ are also generated. This is replicated $20$ times. For each replication, we solve~\eqref{def:maindiscpbprimal} with $g_n$ and the same regularization parameter $\lambda$, and we compute the mean across the 20 replications of the squared-error $\norm{\In-I_N u_N^\star}_{L^2(\O)}^2$. The result is reported in Figure~\ref{fig:signalresultp1}(d). The gray-shaded area indicates one standard deviation of the error over the 20 replications. Again, the average error decreases in agreement with the rate $O\pa{(\log(n)/n)^{1/2}}$ predicted by Theorem~\ref{theo:randomnodesvarational}.

\paragraph{Acknowledgements}{This work was supported by the ANR grant GRAPHSIP. JF was partly supported by Institut Universitaire de France. JF and AE would like to acknowledge support within the EU grant No. 777826, the NoMADs project.}

\bibliographystyle{abbrv}
\bibliography{bibliovariational}
\end{document}